\documentclass[12pt, a4paper]{article}
\usepackage[left=1in,right=1in,top=1.0in,bottom=1in]{geometry}

\usepackage{amsmath,amssymb,amsthm,amsfonts,afterpage,hyperref}
\usepackage{amscd,subeqnarray}

\usepackage{graphicx}
\usepackage{caption}
\usepackage{subcaption}
\usepackage{authblk}
\usepackage{tikz}

\usepackage{array}
\usepackage{wrapfig}
\usepackage{color}
\usepackage{ragged2e}
\usepackage{stmaryrd}
\usepackage{siunitx}
\usepackage{multirow}
\usepackage{xcolor,cancel}
\usepackage{boldfonts}
\usepackage{symbols}
\usepackage{footmisc}
\usepackage{float}
\usepackage{setspace}
\usepackage{bm}
\usepackage{booktabs}
\usepackage{tikz, xcolor}
\usepackage{soul}
\usepackage{pdfpages} 

\usepackage[T1]{fontenc}
\usepackage[utf8]{inputenc}

\usepackage{cite}

\usepackage[linesnumbered,ruled,vlined]{algorithm2e}
\SetKwInput{KwInput}{Input}                
\SetKwInput{KwOutput}{Output} 

\usepackage[numbers,sort&compress]{natbib}
\setcitestyle{square}

\usetikzlibrary{shapes,arrows, positioning,calc}	
\usetikzlibrary{arrows,snakes,backgrounds}

\newcommand{\bfa}[1]{\boldsymbol{#1}} 			
			
\newcommand{\bfeps}{\boldsymbol{\epsilon}}

\newcommand{\Sym}{\text{Sym}}   			%
\newcommand{\ddiv}{\text{div}}     				%
\newcommand{\tr}{\text{tr}}       				%

\DeclareMathAlphabet{\mathpzc}{OT1}{pzc}{m}{it}

\newcommand{\bfn}{\boldsymbol{n}}	
	
\newcommand{\bfu}{\boldsymbol{u}}	
\newcommand{\bfb}{\boldsymbol{b}}	
\newcommand{\bfv}{\boldsymbol{v}}	

\newcommand{\bfF}{\boldsymbol{F}}	
\newcommand{\bfE}{\boldsymbol{E}}	
\newcommand{\bfC}{\boldsymbol{C}}
\newcommand{\bfB}{\boldsymbol{B}}	
\newcommand{\bfx}{\boldsymbol{x}}	
\newcommand{\bfX}{\boldsymbol{X}}	
	
\newcommand{\bfT}{\boldsymbol{T}}		
\newcommand{\bfI}{\boldsymbol{I}}	 
\newcommand{\bfzero}{\boldsymbol{0}}

\newcommand{\bff}{\boldsymbol{f}}	
\newcommand{\bfg}{\boldsymbol{g}}	

\newtheorem{remark}{Remark}

\newtheorem*{cf}{Strong formulation}
\newtheorem*{cwf}{Continuous weak formulation}

\newtheorem*{dfep}{Discrete Finite Element Problem}

\newtheorem*{pia}{Picard's Iterative Algorithm}

\providecommand{\keywords}[1]
{
  \small	
  \textbf{\textit{Keywords---}} #1
}

\title{{A computational model for crack-tip fields in a 3-D porous elastic solid with material moduli dependent on density}}

\author[1]{Kun Gou\thanks{kun.gou@tamusa.edu}}
\author[2]{S. M. Mallikarjunaiah\thanks{m.muddamallappa@tamucc.edu}\thanks{corresponding author}}
\affil[1]{Department of Computational, Engineering, and Mathematical Sciences,
Texas A\&M University-San Antonio,
San Antonio, TX 78224, USA}
\affil[2]{Department of Mathematics \& Statistics,
Texas A\&M University-Corpus Christi,
Corpus Christi, TX 78412, USA}
\date{}

\begin{document}

\maketitle
	    
\begin{abstract}
A mathematical model for crack-tip fields is proposed in this paper for the response of a three-dimensional (3-D) porous elastic solid whose material moduli are dependent on the density. Such a description wherein the generalized Lam\`e coefficients are nonlinear functions of material stiffness is more realistic because most engineering materials are porous, and their material properties depend on porosity and density. The governing boundary value problem for the static equilibrium state in a 3-D, homogeneous, isotropic material is obtained as a second-order, quasilinear partial-differential-equation system with a classical traction-free crack-surface boundary condition. The numerical solution is obtained from a continuous trilinear Galerkin-type finite element discretization. A Picard-type linearization is utilized to handle the nonlinearities in the discrete problem. The proposed model can describe the state of stress and strain in various materials, including recovering the classical singularities in the linearized model. The role of \textit{tensile stress},  \textit{stress intensity factor} (SIF), and \textit{strain energy density} are examined. The results indicate that the maximum values of all these quantities occur directly before the crack-tip, consistent with the observation made in the canonical problem for the linearized elastic fracture mechanics. One can use the same classical local fracture criterion, like the maximum of SIF,  to study crack tips' quasi-static and dynamic evolution within the framework described in this article. 
 \end{abstract}

\noindent \keywords{Implicit constitutive relation, Porous solid, Density-dependent moduli, Finite element method, Crack-tip fields, Tensile stress, Stress intensity factor, Strain energy density   }

\section{Introduction}\label{sec:Intro}

In the classical theory of elasticity, most engineering materials are modeled as compressible. Materials with \textit{Poisson's ratio} equal to $0.5$ do not endure volume changes under any applied deformation. Such materials have no density change; therefore, they must be approximated as incompressible with their constant density. Yet there are many compressible materials, such as polymers, in which density affects the material properties such as strength, hardness, and brittleness. Therefore, one needs to take material moduli dependence on the density to describe the elasticity of polymer networks and the response of such materials to mechanical loading \cite{rubinstein2002elasticity}.  A correction to the classical elasticity model is needed to model the bulk behavior of many porous materials. In the literature, the response of many polymeric materials undergoing infinitesimal deformation is modeled via the classical linearized elastic constitutive relationship in which Young's modulus, Poisson's ratio, and Lam\`e constants are considered constants. In a study on the strength of silica gels \cite{leventis2002nanoengineering}, a nonlinear response is reported between load and strain (see Figure~3 in \cite{leventis2002nanoengineering}) for a three-point bending of aerogels and Young's modulus is computed accordingly. Other studies also demonstrate that the material moduli depend on the density \cite{chandrasekaran2017carbon,anjos2010tensile}. In most solids, the porosity at the macro-level depends on the local density, and the porosity and density are inversely proportional to each other\cite{ALAGAPPAN2023100162}. 

In a series of articles \cite{rajagopal2003,rajagopal2007elasticity,rajagopal2014nonlinear,rajagopal2007response,rajagopal2009class},  Rajagopal generalized the framework of elastic bodies and showed that, through a rigorous linearization, one can obtain constitutive relations that are implicitly related between the linearized strain and Cauchy stress. Such structures are ideal for describing the response of the elastic bodies incapable of converting mechanical work into thermal energy. More importantly, these generalizations are shown to include both Cauchy and Green elastic bodies as a particular subclass of elastic material models \cite{bridges2015implicit,Mallikarjunaiah2015,MalliPhD2015,ortiz2014numerical,ortiz2012}.   The aforementioned nonlinear behavior can be modeled within the implicit theory of elasticity \cite{rajagopal2014nonlinear,muliana2018determining,kowalczyk2019finite}.  Rocks \cite{bustamante2020novel} and many rubber-like materials \cite{bustamante2021new} can be modeled using a particular subclass of models within the larger class of nonlinear relationships from \cite{rajagopal2003,rajagopal2007elasticity}. In \cite{bustamante2021new},  a constitutive model for rubber is proposed using the principal stress as the main variables, wherein a better corroboration with experimental data can be obtained than that for Ogden's model \cite{ogden1972large}.

The nonlinear relations proposed in \cite{rajagopal2003,rajagopal2007elasticity} between the strain and stress cannot be obtained from Cauchy or Green elasticity by applying the same standard assumption that the displacement gradients are infinitesimal. From the experimental viewpoint, the classical linear elasticity can't adequately explain the nonlinear behavior between the stress and strain reported in \cite{saito2003multifunctional}. During the cold swagging process,  for the metal alloys Ti-12Ta-9Nb-3V-6Zr-O and 23Nb-0.7Ta-2Zr-O, the nonlinear behavior is observed well within a 2.5\% order of strain. A similar nonlinear behavior between stress and strain has been reported in other experimental studies on titanium alloys \cite{withey2008deformation}, gum metals \cite{zhang2009fatigue}, and rocks \cite{johnson1996manifestation}.  Several materials such as gum metal \cite{kulvait2019state}, titanium alloys (e.g., the one used in the orthopedic applications \cite{tian2015nonlinear,devendiran2017thermodynamically}), and other alloys (e.g., Ti-30Nb-10Ta-5Zr (TNTZ-30)) are known to exhibit an apparent nonlinear behavior within about $1\%$ strain \cite{hao2005super,saito2003multifunctional}. Such behavior is traditionally modeled using a linear relationship between the infinitesimal strain and Cauchy stress under the classical linearized elasticity construct, which is inadequate because a linear relationship can't fit the experimentally observed nonlinear relations.

The new theory of elasticity proposed by Rajagopal  \cite{rajagopal2003,rajagopal2007elasticity} has been used to revisit the classical problems in fracture mechanics, such as plane-strain and anti-plane static cracks in \cite{rajagopal2011modeling,gou2015modeling,HCY_SMM_MMS2022}, plane-strain cracks in a transversely isotropic elastic body \cite{Mallikarjunaiah2015}, quasi-static crack propagation in a new class of material models \cite{yoon2021quasi,lee2022finite}, static cracks in the implicit thermo-elastic body \cite{yoon2022finite}, and static cracks in 3-D materials body \cite{gou2023MMS,gou2023computational}. One of the significant advantages of the nonlinear constitutive theory in \cite{rajagopal2003,rajagopal2007elasticity} to study cracks and fractures is that the strain values are consistent with these from the linearization procedure used to derive the relationship, contrary to the classical square-root singularity predicated by the canonical boundary value problem (BVP) within the linear elasticity. However, one can still find the stress concentration within the framework proposed in \cite{rajagopal2003,rajagopal2007elasticity} suggesting using the same local fracture criterion as in the linear elastic fracture mechanics.

A subclass of nonlinear models in \cite{rajagopal2003,rajagopal2007elasticity} for the response of elastic materials such as rocks and concrete wherein the material moduli dependent upon the density are obtained in \cite{rajagopal2021b,rajagopal2022implicit}. In such models, both linearized strain and Cauchy stress appear linearly \cite{itou2022investigation,itou2021implicit,itou2023generalization,pruuvsa2022pure,pruuvsa2023mechanical}. The material response of the models in \cite{rajagopal2021b,rajagopal2022implicit} is not linear as the constitutive relations contain the bilateral product of Cauchy stress and linearized strain. However, the underlying structure of such models allows one to describe the linear homogeneous material with the inherent properties depending upon the density. In \cite{HYCSMM2023}, a two-dimensional model for materials with properties dependent upon the density is proposed, and the crack-tip fields along with mechanical stiffness for the anisotropic material are studied. A variety of phenomena were captured for both tensile and in-plane shear loading. Several engaging failure scenarios were reported, which were impossible to obtain under the classical linear elasticity models. 

In this paper, we formulate a model for a 3-D plate with an edge crack, i.e., a v-notch with a slight opening angle. A BVP is developed for the v-shaped crack geometry within the setting of nonlinear constitutive relations with both Cauchy stress and strain appearing linearly. A convergent numerical method based on the continuous Galerkin-type finite elements is proposed for the quasi-linear elliptic partial differential equation system. Several postprocessing items, including the crack-tip stress and strain,  stress intensity factor, and strain energy density, are presented along a line leading to the v-notch tip. Our study is a basis for the future development of fracture theory for brittle materials with density-dependent moduli that play a crucial role in incorporating the material response to mechanical and thermal stimuli.

\section{Constitutive modeling}
\subsection{Basic notations}
Here, we briefly present the geometrical description of the elastic body under investigation. We assume that the body contains bulk material around a crack tip with the ultimate goal of resolving the crack-tip fields.  Let ${\Omega}$ be a region in the 3-D Euclidean space $\mathbb{R}^3$ occupied by the body as the reference configuration. The Lipschitz continuous boundary of $\Omega$ is decomposed into a disjoint union of a Neumann boundary $\Gamma_N$ and a nonempty Dirichlet boundary $\Gamma_D$ (the boundary measure of $\Gamma_D >0$), i.e. $\partial \Omega = \overline{\Gamma_N} \cup \overline{\Gamma_D}$. Also, let $\bfa{n} = \left( n_1, \; n_2, \; n_3 \right)$ be the outward unit normal. We suppose that $\Gamma_c$ represents the {V-shaped} surface completely inside $\overline{\Omega}$. Let $\bfX = \left(X_1, \,  X_2, \, X_3 \right)$  and  $\bfx = \left(x_1, \, x_2, \, x_3 \right)$ denote typical points in the reference and deformed configurations of the body, respectively.  Let $\bfu \colon \Omega \to \mathbb{R}^3$ denote the displacement field satisfying
 \begin{equation}
 \bfu = \bfx - \bfX.  
 \end{equation} 
Let $\bfx = \chi(\bfX, \; t)$ be the smooth deformation of the body ${\Omega}$, where $t$ denotes the time. The displacement gradients $\dfrac{\partial \bfu}{\partial \bfX}$ and $\dfrac{\partial \bfu}{\partial \bfx}$ are given by
\begin{equation}
\dfrac{\partial \bfu}{\partial \bfX} = \nabla_{\bfX} \bfu = \bfa{F} - \bfa{I}, \quad \mbox{and} \quad \dfrac{\partial \bfu}{\partial \bfx} = \nabla_{\bfx} \bfu =  \bfa{I} - \bfa{F}^{-1},
\end{equation}
where $\bfF \colon \Omega \to \mathbb{R}^{3 \times 3}$ denotes the deformation gradient.  

Let $\Sym(\mathbb{R}^{3 \times 3})$ be the space of $3 \times 3$ symmetric tensors equipped with the inner product $\bfa{M} \colon \bfa{N} = \sum_{i, \, j=1}^3 \, \bfa{M}_{ij} \, \bfa{N}_{ij}$ for all $\bfa{M} $ and $\bfa{N} $ in $\Sym(\mathbb{R}^{3 \times 3})$, and with  {the associated norm}  {$\| \bfa{M} \| = \sqrt{\bfa{M} \colon \bfa{M}}$}. 

The \textit{right Cauchy-Green stretch tensor} $\bfC \colon \Omega \to \mathbb{R}^{3 \times 3}$, the \textit{left Cauchy-Green stretch tensor} $\bfB \colon \Omega \to \mathbb{R}^{3 \times 3}$, the \textit{Lagrange strain} $\bfE \colon \Omega \to \mathbb{R}^{3 \times 3}$, and the \textit{ linearized strain tensor} $\bfeps \colon \Omega \to  \Sym(\mathbb{R}^{3 \times 3})$ are defined by
\begin{subequations}
\begin{align}
 \bfB &:=\bfF\bfF^{\mathrm{T}}, \; \bfC :=\bfF^{\mathrm{T}}\bfF,  \; \bfE := \dfrac{1}{2} \left( \bfC - \bfI \right),\\
\bfeps &:= \dfrac{1}{2} \left( \nabla \bfu + \nabla \bfu^{\mathrm{T}}\right),
\end{align}
\end{subequations}d
where $\left( \cdot \right)^{\mathrm{T}}$ denotes the \textit{transpose} operator for tensors, and $\bfI$ is the
3-D identity tensor. Our interest in this paper is to derive the mathematical models for the response of porous elastic solids undergoing only small displacement gradients, and thus, it is required that 
\begin{equation}\label{small_grad}
\max_{\bfX \in \Omega} \| \nabla \bfu \|  \ll \mathcal{O}(\delta), \quad \delta \ll 1.
\end{equation}
Under the above assumption \eqref{small_grad}, one gets 
\begin{subequations}\label{lin_results}
\begin{align}
&\bfB \approx \bfI + 2 \bfeps + \mathcal{O}(\delta^2), \;\;  \bfC \approx \bfI + 2 \bfeps + \mathcal{O}(\delta^2), \;\;  \bfE \approx  \bfeps,  \\
&\det \bfF = 1 + \tr(\bfeps). \label{detF_linear}
\end{align}
\end{subequations}
Let $\bfT \colon \Omega \to  \Sym(\mathbb{R}^{3 \times 3})$ be the \textit{Cauchy stress tensor} in the current configuration and for the most general case of an isotropic homogeneous compressible elastic solid, the Cauchy stress takes the following form
\begin{equation}\label{Tequation}
\bfT = \widehat{\beta_0} \bfI + \widehat{\beta_1} \bfB +  \widehat{\beta_2} \bfB^2,
\end{equation}
where the material moduli  $\widehat{\beta_i}, \; i=0, \; 1, \; 2$ are functions of $\tr \, \bfB$, $\quad \dfrac{1}{2} \left[ \left(   \tr \, \bfB   \right)^2 - \tr \, \bfB^2  \right]$ and $\det \bfB$. The quantity $\widehat{\beta_0}$ can be the spherical stress due to the constraint of incompressibility.  The Cauchy stress $\bfT$ satisfies the linear momentum balance and the angular momentum balance, respectively, as below
\begin{subequations}\label{balance_equations}
\begin{align} 
\rho \, {\ddot{\bfu}} &= \ddiv \, \bfT + \rho \, \bfb, \\
\bfT &= \bfT^{\mathrm{T}}.
\end{align}
\end{subequations}
Here, $\rho$ is the density of the material in the current configuration, and $\bfb \colon \Omega \to \mathbb{R}^3$ is the body force density.

In this paper, we only consider the case of infinitesimal displacement gradients. Hence, we are well within the realm of linearized elasticity. However, we will work with algebraically nonlinear constitutive relationships to describe the response of geometrically linear material.  In the following subsection, we briefly describe the derivation of nonlinear relations between linearized strain and Cauchy stress, but the material parameters depend on density.

\subsection{Implicit constitutive relations}\label{implirela}
In this study, we aim to characterize the behavior of 3-D porous solids whose material moduli depend on the density. To this end, we closely follow the work of Rajagopal in \cite{rajagopal2003,rajagopal2007elasticity,rajagopal2021a,rajagopal2021b} where an implicit relation among the stress $\bfT$, deformation gradient $\bfF$, and density $\rho$ is used to describe the response of solids, satisfying
\begin{equation}\label{eq:implicit_relation}
\bff(\rho, \, \bfT,\,\bfB) =\bf0,
\end{equation}
where the tensor-valued function $\bff$ is assumed to be \textit{isotropic}. The above relationship in its most general representation form can  be written as \cite{gokulnath2017representations}
\begin{align} 
0 &= \beta_0 \, \bfI + \beta_1 \, \bfT + \beta_2 \, \bfB +  \beta_3 \, \bfT^2 + \beta_4 \, \bfB^2 +  \beta_5 \, ( \bfT \bfB + \bfB \bfT) + \beta_6 \, ( \bfT^2 \bfB + \bfB \bfT^2)  \notag \\
& + \beta_7 \, ( \bfB^2 \bfT + \bfT \bfB^2) +  \beta_8 \, ( \bfT^2 \bfB^2 + \bfB^2 \bfT^2),  
\end{align}
where the scalar-valued functions $\beta_i$, for $i=0,\, 1, \,  \ldots, 8$, are material moduli dependent on the density and  invariants of $\bfT$ and $\bfB$ as below
\begin{align*}
\big\{ &\rho, \, \tr (\bfT), \, \tr(\bfB), \, \tr(\bfT^2), \, \tr(\bfB^2), \,\tr(\bfT^3), \, \tr(\bfB^3), \,\tr(\bfT \, \bfB), \, \tr(\bfT^2 \, \bfB) \notag \\
&\tr(\bfT \, \bfB^2), \, \tr(\bfT^2 \, \bfB^2)   \big\}.
\end{align*}
A special sub-class of the above relations is given by \cite{murru2021stress,murru2021ZAMM,vajipeyajula2023stress,itou2021implicit,itou2022investigation}  
\begin{equation}\label{imp_model}
\widetilde{\delta}_0 \, \bfI + \widetilde{\delta}_1 \, \bfeps + \widetilde{\delta}_2 \, \bfT + \widetilde{\delta}_3 \, \bfT^2 + \widetilde{\delta}_4 \, \left( \bfeps \, \bfT + \bfT \, \bfeps \right) + \widetilde{\delta}_5 \, \left( \bfeps \, \bfT^2 + \bfT^2 \,\bfeps \right)    =0,
\end{equation}
where the functions  $\widetilde{\delta}_i, \; i=0, \, 2, \, 3$ depend linearly upon the invariants of $\bfeps$ and arbitrarily upon the invariants of  $\bfT$. While, the functions  $\widetilde{\delta}_i, \; i=1, \, 4, \, 5$ depends on the invariants of $\bfT$. The above sub-class of relations have been studied extensively in the literature for a variety of BVPs \cite{HYCSMM2023,yoon2022finite,HCY_SMM_MMS2022,yoon2021quasi,lee2022finite,ALAGAPPAN2023100162,murru2022density,gou2015modeling}.  

The density is considered globally constant in all the studies concerning the development of mathematical models for the response of solid bodies containing cracks and fractures. The balance of mass (or continuity equation) in the material description yields 
\begin{equation}\label{mass_balance1}
\rho_0 = \rho \, \det \bfF,
\end{equation}
where $\rho_0$ is the reference density. By \eqref{detF_linear},  \eqref{mass_balance1} reduces to 
\begin{equation}\label{mass_balance2}
\rho_0 = \rho \, \left( 1 + \tr(\bfeps)  \right).
\end{equation}
An implicit constitutive relation describing the small strain response of the porous elastic solid  is given by a special subclass of \eqref{imp_model}, \cite{itou2021implicit,HYCSMM2023}
\begin{equation} \label{spe_model1}
 \bfeps = C_1 \, \left( 1 + \xi_1 \, \tr \, \bfeps  \right) \; \bfT + C_2 \, \left( 1 + \xi_2 \, \tr \, \bfeps  \right) \;  \left(  \tr \bfT \right)\; \bfI, 
\end{equation}
where $\xi_1, \, \xi_2, \, C_1, \, C_2$ are all constants. {Notice that the classical linearized elasticity model is a special case of the above relation \eqref{spe_model1} for $\xi_1\to 0$ and $\xi_2\to 0$. Hence, the model \eqref{spe_model1} can generate broader classes of constitutive relations for different values of the parameters.} The constants $C_1$ and $C_2$ are shown to be
\begin{equation}\label{eq:material_coeff}
C_1 = \dfrac{1 + \nu}{E} = \dfrac{1}{2 \, \mu} >0, \quad C_2 = - \dfrac{ \nu}{E}  <0,
\end{equation}
where $E$ is the Young's modulus, $\mu$ is the shear modulus (one of the two Lam\`e constants), and $\nu$ is the Poisson's ratio. For linearized elasticity, the Lam\`e constants $\lambda$ and $\mu$ are related to Young's modulus  and  Poisson's ratio via
\begin{equation}\label{eq:linear_lame}
\lambda =  \dfrac{ E \, \nu}{(1 + \nu) (1 - 2 \nu)}, \quad \mu = \dfrac{ E }{2 (1 + \nu)}.
\end{equation}
The above constitutive class (\ref{spe_model1}) provides a relationship between the linearized strain and the Cauchy stress. 

\begin{remark}
The expression $\left( 1 + \xi_1 \, \tr \, \bfeps  \right)$ can be recognized as density-dependent moduli since  $\tr \, \bfeps$ can be expressed in terms of densities of two configurations. Thus, model \eqref{spe_model1} can describe the response of porous materials in which the density changes predominantly, and also, the material moduli of such solids generally depend on the density. On the contrary, the material moduli are density-independent under the linearized elastic description of solids. 
\end{remark}

\begin{remark}
A large body of literature is devoted to the BVP formulated within the realm of the linearized theory of elasticity \eqref{Tequation}-\eqref{balance_equations} for the brittle elastic bodies containing cracks and fractures. Based on the choice of the invariants in  \eqref{Tequation} and under the valid assumption of linearized elasticity \eqref{small_grad}, the stress and strain may be related linearly to study quasi-static evolution of crack-tip \cite{knees2015quasilinear,negri2017quasi,kuttler2006quasistatic}.   As in problems of cracks with classical traction-free crack-surface boundary conditions, the strain becomes very large in the neighborhood of crack-tip, thereby contradicting the fundamental doctrine under which the theory is built \cite{rajagopal2014nonlinear}. However, there have been a lot of modifications proposed to the canonical BVP within the context of a linearized theory \cite{sendova2010new,walton2012note,ferguson2015numerical,sinclair2004stress,kim2011analysis}. The implicit theory of elasticity proposed by Rajagopal  \cite{rajagopal2003} attempts to build in new algebraic nonlinear relations for the elastic solids, and the BVPs within this new paradigm have been shown to predict bounded strain in the neighborhood of the crack-tip \cite{rajagopal2011modeling,gou2015modeling}. 
\end{remark}

In the next section, using the novel constitutive relations, we develop a mathematical model to describe the state of stress and strain in porous elastic materials. Subsequently, we present a stable discretization of the partial differential equation systems and demonstrate the difference between the numerical results from the linearized elasticity model and the novel implicit model. 

\section{BVP for 3-D porous elastic solid}\label{sec:BVP}
In this section, we develop a BVP for the response of a 3-D porous elastic solid described by the special constitutive relation in Sec. \ref{implirela}. The main objective is to resolve the crack-tip fields within the new relationship wherein the stress and linearized strain appear linearly, and the material moduli of the porous material depend on the density. We assume the material bulk body under investigation is homogeneous, isotropic, and initially unstrained and unstressed.  To describe the state of stress-strain in the new class of models, we consider a static problem. The balance of linear momentum in the absence of body force reduces to 
\begin{equation}
 - \nabla \cdot \bfT  = \boldsymbol{0} \quad \text{in} \quad \Omega.  \label{eq_blm}
\end{equation}
To model the stress response of the material whose material moduli are dependent upon the density, we utilize the constitutive relationship described in \eqref{spe_model1} under 
\[
\xi_1=\xi_2=\beta,  \quad \beta \in \mathbb{R}. 
\]
Such a constitutive relationship is invertible to express the Cauchy stress as a function of the strain as
\begin{equation}\label{def-T}
 \bfT=\frac{\mathbb{E}[\bfeps]}{1 + \beta \, \tr(\bfeps)},
\end{equation}
 where the fourth-order elasticity tensor $\mathbb{E}[\bfeps]$ satisfies
\begin{equation}
\mathbb{E}[\bfeps] := \dfrac{E}{(1+ \nu)} \bfeps + \dfrac{ \nu \, E}{(1 -2\nu)  (1 + \nu) } \tr \bfeps \; \bfI.
\end{equation}
In the new description of Cauchy stress for porous elastic materials, the generalized Lam\`e coefficients for the nonlinear constitutive model are
\begin{equation}\label{eq:c1_c2_dd_model}
\lambda := \dfrac{\nu \, E}{(1 + \nu)(1-2\nu)(1+ \beta \, \tr(\bfeps))}, \quad \mu := \dfrac{E}{2(1 + \nu)(1+ \beta \, \tr(\bfeps))}.
\end{equation}
Notice that the stress function defined in \eqref{def-T} is nonlinearly dependent upon the strain due to the function itself and the nonlinearity of the Lam\`e coefficients. 
\begin{remark}
It is clear from the constitutive relation \eqref{spe_model1} that when $\bfT =0$, $\bfeps=0$, and vice versa,  both as expected. Furthermore,  $ \| \bfT \| \to +\infty$ implies $\tr \, \bfeps \to -1/\beta$, generating a limited strain behavior. The strain limit needs to be positive for tensile stress; hence, for stress concentration problems such as cracks and fractures under tension, the constant $\beta$ should be negative.  Conversely, for the compressive stresses, the above model shares the same ``negative'' effects as the classical linearized model regarding the strain singularity near the crack tips. Hence, one must ensure that the linearized strain is appropriately small for the compressive stress problems. Overall, the model proposed in  \eqref{spe_model1} is the most generalized one, and more importantly, the strain is limited in tension. Still, the stress has to be appropriately small so that the strain is infinitesimal in compression. 
\end{remark}

Utilizing  \eqref{spe_model1},  \eqref{eq_blm}, and \eqref{def-T}  we obtain the following BVP:
\begin{cf}
Given all the material parameters, find $\bfu = (u_1, \, u_2, \, u_3)$ such that 
\begin{subequations}\label{f1}
\begin{align}
- & \, \nabla \cdot \left[    \frac{\overline{c}_1/2  \, \left( \nabla \bfu + \nabla \bfu^T \right) + \overline{c}_2\nu \, \nabla \cdot \bfu  \, \bfI}{1 + \beta \;  \nabla \cdot \bfu}         \right]  = \bff, \; \mbox{in} \;\; \Omega, \;\;  \label{f1_1} \\
 \bfeps &= \dfrac{(1+\nu)(1+ \beta \, \tr(\bfeps))}{E} \, \bfT - \dfrac{\nu(1+ \beta \, \tr(\bfeps))}{E} \, \tr(\bfT) \,  \bfI  \;\; \text{in} \;\; \Omega,  \\
 \bfu &= \widehat{\bfu}, \quad  \mbox{on} \quad \Gamma_D \\
 \bfT  \bfn &= \bfg, \;\;  \mbox{on} \;\; \Gamma_N,
\end{align}
\end{subequations}
where  $\overline{c}_1$ and $\overline{c}_2$ are given by
\begin{equation}
\overline{c}_1 = \dfrac{E}{(1+\nu)}, \quad \mbox{and} \quad \overline{c}_2= \dfrac{ E}{(1+\nu)(1 - 2 \nu)},
\end{equation}
$\bfg \colon \Gamma_N \to \mathbb{R}^3$ is the given traction, and $\widetilde{\bfu} \colon \Gamma_D \to \mathbb{R}^3$ is the given boundary displacement.
\end{cf}

The above BVP  \eqref{f1} is nonlinear; currently, no analytical tools are available to obtain its closed-form solution. Therefore, we resort to a convergent numerical method, such as the finite element method, to find the approximate solution.  In the next section, we propose a stable finite-element discretization for (\ref{f1}) coupled with a Picard-type iteration algorithm to obtain the numerical solution.

\section{Finite element discretization}
We aim to study the response of the 3-D porous elastic solid whose material moduli depend upon the density. The constitutive relationship appears linear in the stress and linearized strain \cite{murru2021stress,itou2021implicit,rajagopal2021b,HYCSMM2023}.  Suppose that the domain $\Omega$ is bounded and connected. Let $L^p(\Omega)$ and $W^{k,p}(\Omega)$ denote the standard Lebesgue and Sobolev spaces in 1-D, respectively.  The corresponding spaces for the 3-D vector-valued functions and symmetric $3\times3$ tensor-valued functions are denoted by  $(L^p(\Omega))^3$ and $(L^p(\Omega))^{3\times3}_{\Sym}$ for the Lebesgue space, and $(W^{k,p}(\Omega))^3$ and $(W^{k,p}(\Omega))^{3\times3}_{\Sym}$ for the Sobolev space. To approximate the displacement that vanishes on the boundary $\Gamma_D$, we consider the Soboleve space $(W^{1,p}_0(\Omega))^3$, where $W_0^{1,p} $  means the closure of the set of infinitely differentiable functions with compact support in $\Omega$ and its norm defined as for the space    $W^{1,2}(\Omega)$, satisfying
\begin{equation}
W_0^{1,p} := \overline{C_{c}^{\infty}(\Omega)}^{ \| \cdot \|_{1,p}} \label{def-H01}.
\end{equation}
We also define the following subspaces of $\left( W^{1,\,2}(\Omega)\right)^3$:
\begin{subequations}
\begin{align}
V_{\bfzero} &:= \left\{ \bfu \in \left( W^{1,\,2}(\Omega)\right)^3 \colon \; \bfu=\bfzero \quad  \mbox{on} \;\;\Gamma_D\right\}, \label{test_V0}\\
V_{\widehat{\bfu}} &:= \left\{ \bfu \in \left( W^{1,\,2}(\Omega)\right)^3 \colon \; \bfu=\widehat{\bfu} \quad \mbox{on} \;\; \Gamma_D\right\}. \label{test_Vu0}
\end{align}
\end{subequations}
To obtain the weak formulation of the proposed BVP, we multiply  \eqref{f1_1} by a test function $\bfv \in {V}_{\bfzero}$ and integrate by parts. Then, using Green's formula along with appropriate boundary conditions, we obtain the following weak formulation. 
\begin{cwf}
Find  $\bfu \in {V}_{\widehat{\bfu}}$, such that 
\begin{equation}\label{eq:weak_formulation}
    a(\bfu, \, \bfv) = l(\bfv), \;  \forall\, \bfv \in {V}_{\bfzero},
\end{equation}
where the bilinear term $a(\bfu, \, \bfv)$ and the linear term $l(\bfv)$ are defined by
\begin{subequations}\label{def:A-L}
\begin{align}
    a(\bfu, \, \bfv) &= \int_{\Omega} \left[    \frac{\overline{c}_1/2  \, \left( \nabla \bfu + \nabla \bfu^T \right) + \overline{c}_2\nu \, \nabla \cdot \bfu  \, \bfI}{1 + \beta \;  \nabla \cdot \bfu}         \right] \colon \bfeps( {\bfv}) \;   d\bfx\, ,  \label{eq42a}\\
 l (\bfv) &= \int_{\Omega} \bff \cdot \, \bfv \; d\bfx  +  \int_{\Gamma_N} \bfa{g} \cdot \bfa{v} \; ds.  
\end{align} 
\end{subequations} 
\end{cwf}
The above continuous weak formulation is nonlinear, and hence, obtaining a sequence of linear problems for numerical simulations to facilitate investigating the stress-strain concentration near the crack tip in a 3-D body would be helpful. We construct a \textit{Picard's type iterative algorithm} at the continuous level and then discretize the subsequent linear problems using the stable Galerkin-type finite-element method by the following algorithm.
\begin{pia}
For $n=0, \, 1, \, 2, \, \ldots$, and given $\bfu^{0} \in {V}_{\widehat{\bfu}}$,  find  $\bfu^{n+1} \in {V}_{\widehat{\bfu}}$, such that 
\begin{equation}\label{eq_pia}
    a(\bfu^{n}; \; \bfu^{n+1}, \, \bfv) = l(\bfv), \;  \forall\, \bfv \in {V}_{\bfzero},
\end{equation}
where
\begin{equation}\label{eq42a}
    a(\bfu^{n}; \; \bfu^{n+1}, \, \bfv) = \int_{\Omega} \left[    \frac{\overline{c}_1/2  \, \left( \nabla \bfu^{n+1} + \nabla (\bfu^{n+1})^T \right) + \overline{c}_2\nu \, \nabla \cdot \bfu^{n+1}  \, \bfI}{1 + \beta \;  \nabla \cdot \bfu^{n}}         \right] \colon \bfeps( {\bfv}) \;   d\bfx.
\end{equation} 
\end{pia}



\subsection{Finite element space and discrete problem}
We briefly describe the finite element discrete problem corresponding to the defined continuous weak formulation \eqref{eq_pia}-\eqref{eq42a}. We assume that the 3-D material domain $\Omega$ is a non-convex polyhedral domain and that there exists a mesh $\mathcal{T}_h$ with $h>0$ being the mesh size and such a mesh being either quasi-uniform or a priori refined mesh. The finite element discretization of the domain $\Omega$ is assumed to be 
conforming and shape-regular in the sense of  Ciarlet \cite{ciarlet2002finite}. The discretization is done in such a way that for any two elements $K_1, K_2 \in \mathcal{T}_h$, $\overline{K}_1 \cap \overline{K}_2$ can be one of several cases including a null set, a vertex, an edge, and even $\overline{K}_1$ (or $\overline{K}_2$),  and $\bigcup\limits_{K \in \mathcal{T}_h} \overline{K} = \overline{\Omega}$. 

 For approximating the  displacement field  $\bfu$,  we define the following space,
\begin{equation}
S_h = \left\{  \bfu_h \in  \left( C(\overline{\Omega})\right)^3 \colon \left. \bfu_h\right|_K \in \mathbb{Q}_k^d, \; \forall K \in \mathcal{T}_h \right\},
\end{equation}
where $\mathbb{Q}_k^d$ is a set of tensor-products of polynomials up to an order of $k$ over the reference cell $\widehat{K}$. Then the discrete  approximation space is
\begin{equation}\label{app-spaces}
\widehat{V}_h = S_h \, \cap  \, V_{\widehat{\bfu}}. 
\end{equation}
The discrete finite element problem is given below.
\begin{dfep}
Given $\beta$, the Dirichlet boundary data $\bfu^{0}_h \in \widehat{V}_h$, and the $n^{th}$ iteration solution $\bfu^n_h \in \widehat{V}_h$, for $n=0, 1, 2, \cdots $, find $\bfu^{n+1}_h  \in \widehat{V}_h$ such that 
\begin{equation}\label{discrete-wf}
   a(\bfu_h^n; \, \bfu^{n+1}_h,\, \bfv_h) = l(\bfv_h), \forall\, \bfv_h \in \widehat{V}_h,  
\end{equation}
where the bilinear and linear terms are given by
\begin{subequations}\label{A-L-Def}
\begin{align}
a(\bfu_h^n; \, \bfu^{n+1}_h,\, \bfv_h) &=\int_{\Omega} \left[    \frac{\overline{c}_1/2  \, \left( \nabla \bfu_h^{n+1} + \nabla (\bfu_h^{n+1})^T \right) + \overline{c}_2 \nu \, \nabla \cdot \bfu_h^{n+1}  \, \bfI}{1 + \beta \;  \nabla \cdot \bfu_h^{n}}         \right] \colon \bfeps( {\bfv_h}) \;   d\bfx\, ,  \label{disc_A} \\
l(\bfv_h) &= \int_{\Omega} \bff \cdot \, \bfv_h \; d\bfx  +  \int_{\Gamma_N} \bfa{g} \cdot \bfv_h \; ds.  \label{disc_L}  
\end{align}
\end{subequations}
\end{dfep}

\begin{remark}
A proper initial guess $\bfu^{0}_h$ is required for fast convergence of the above iterative algorithm. In our implementation, we solved the linear problem first (obtained from \eqref{discrete-wf} by taking $\beta=0$) and subsequently used the computed solution as an initial guess for Picard's iterations.  For multiple values of $\beta$, all our numerical simulations converged within a reasonable number of iterations.  
\end{remark}

\subsection{Algorithm for computation}
 The following algorithm depicts the overall discrete finite element computational procedure to obtain the numerical solution of the BVP.
 
\begin{algorithm}[H]
\SetAlgoLined
\KwInput{Choose the parameters: $\beta, \, E, \, \nu$, $I_{Max}$ (maximum of the iteration number), $Tol$}
Start with a sufficiently refined mesh; \\
Solve the linear problem with $\beta =0$ in \eqref{discrete-wf} to obtain the initial guess $\bfu^{0}_h \in \widehat{V}_h$ for subsequent Picard's iteration algorithm; \\
\For {$n=0, \,1, \, 2, \, \ldots$}{
 \While{$[\text{Iteration Number} < \text{$I_{Max}$}]$ AND $[\text{Residual} > \text{$Tol$}]$}{
  Assemble Equations~\eqref{disc_A} and \eqref{disc_L} using test functions from the discrete finite element space $S_h$ \;
 Use a \textit{direct solver} to solve for $ \bfu_h^{n+1}$\;
  }
  }  
 Save the final converged solution $\bfu_h$ to output files for post-processing\;
 Compute the crack-tip fields (e.g., stress and strain) for visualization.
 \caption{Finite element algorithm for the nonlinear BVP}
 \label{algo001}
\end{algorithm}
 Two important post-processing items  need to be computed  to characterize the crack-tip fields in the porous elastic solid using the final converged solution $\bfa{u}_h$: 
\begin{subequations}
\begin{align}
\bfa{T}_h &=\frac{ \overline{c}_1/2 \, \left( \nabla \bfa{u}_h + (\nabla \bfa{u}_h)^T \right)  + \overline{c}_2\nu \, \nabla \cdot \bfa{u}_h \, \bfa{I} }{1 + \beta \;  \nabla \cdot \bfu_h}, \\
\bfeps_h &= \dfrac{(1+\nu)(1+ \beta \,  \nabla \cdot \bfu_h)}{E} \, \bfT - \dfrac{\nu(1+ \beta \, \nabla \cdot \bfu_h )}{E} \, \tr(\bfT) \,  \bfa{I}.
\end{align}
\end{subequations}
Another important crack-tip field variable that symbolizes the strength of the material and quantifies the crack-tip evolution is the \textit{strain energy density} which is obtained by taking the inner tensor product of the stress $\bfa{T}_h$ and the strain $\bfeps_h$. The current work is the first attempt to present the entire behavior of the 3-D porous elastic solid whose material properties depend upon the density. 


\section{{Computational results}}\label{sec:NumExp}

\begin{figure}[H]
	\centering
	\includegraphics[width=0.5\textwidth]{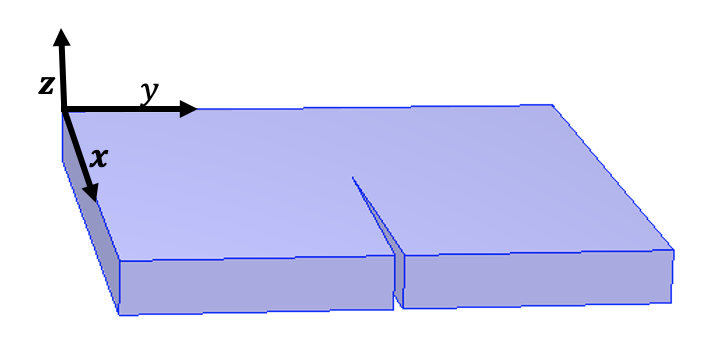}
	\caption{\footnotesize{The geometry of the crack. The crack is on a square plate. The crack angle is $2^o$. The $xyz$ rectangular coordinate system is displayed on the geometry.} }
	\label{plategeo}
\end{figure}

As shown in Fig. \ref{plategeo}, a geometry used in \cite{gou2023MMS} for studying a strain-limited model, we consider a crack embedded in a thin square plate. The crack bears an angle of $2^o$, and the crack tip is on the center of the square plate. Each side of the square is 100 mm long, and the plate is 10 mm thick. The parameter values of the material are $E=10^4$ Pa and $\nu=0.3$. 

\begin{figure}[H]
	\centering
	\subfloat[\footnotesize{Tensile-displacement boundary conditions.}]{\label{tensile}\includegraphics[width=0.5\textwidth]{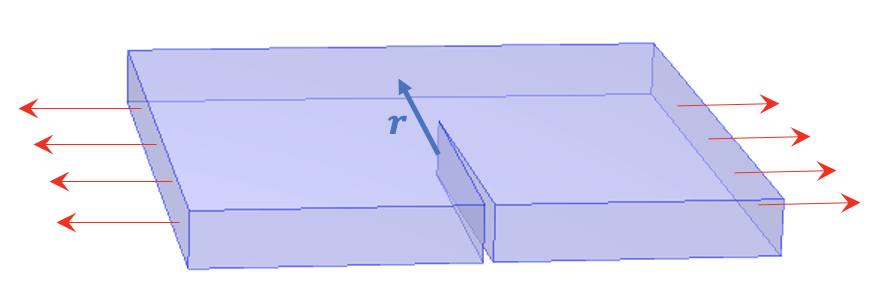}} 
	\subfloat[\footnotesize{Tetrahedral meshes. }]{\label{meshes}\includegraphics[width=0.35\textwidth]{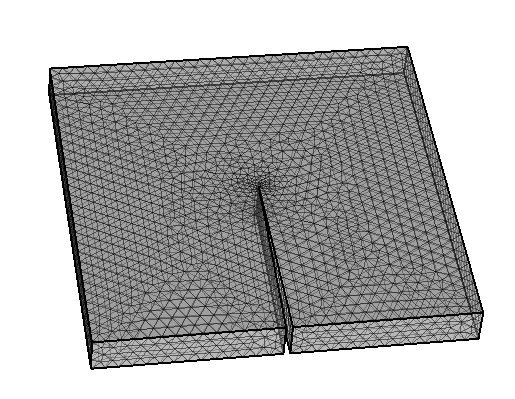}} 
	\caption{{\footnotesize{The tensile displacement boundary conditions are applied on two sides of the plate parallel to the $x$-axis. Other boundaries are traction-free. Tetrahedral meshes are prescribed over the cracked plate. To produce more accurate computational results, finer meshes are created near the crack tip.}}}
	\label{D3Stress22Bound1}
\end{figure}

 Tensile displacement boundary conditions are applied to the plate as shown in Fig. \ref{tensile}. Over the two side boundaries, the displacement magnitude in the $y$-direction is 1 mm. All the other sides are traction-free. We mainly study computational results near the crack tip to understand how the nonlinear model impacts the crack's strain and stress distribution, the stress intensity factor (SIF), and strain energy. In Fig. \ref{tensile}, a blue arrowed line with coordinate $r$ starting from the middle of the edge of the crack tip and perpendicular to the plate's back is created. At the crack tip, $r=0$, and $r$ increases when further away from the crack tip.  Part of the computational results are demonstrated near the crack tip along the line for a pronounced view. Another benefit of choosing the middle line is to minimize the effect of boundary conditions on the computational results. 
 
  In our computations, the \textit{tetrahedral meshes} are prescribed over the whole geometry (Fig. \ref{meshes}). The results near the crack tip are essential to studying the nonlinear model effect. Thus, highly fine meshes are generated in the neighborhood of the crack tip. There are 38052 tetrahedral meshes with 7853 mesh vertices.  Over the boundaries, there are 5718 triangle meshes. The quality of the meshes is measured by skewness, i.e.,  the angular measure of mesh quality compared to the angles of ideal mesh types. The minimum mesh quality is 0.2298, and the average is 0.6754. The shape function of the finite element is quadratic Lagrange. The number of degrees of freedom solved for is 169848.
  The computations are implemented in the Intel(R) Core(TM) i9-9980HK CPU at 2.40 GHz of a Macbook Pro laptop. The physical memory used is 3.69 GB, with a virtual memory of 42.62 GB. We demonstrate the computational results for both positive and negative $\beta$ values at $\beta=-30, \, -20, \,-10, \,0, \,10, \,20, \,30$. {The average computational time for each $\beta$ value is about 38 seconds.} $\beta=0$ generates results for the linear-elasticity fracture model, and the results from $\beta=0$ are compared with results for other  $\beta$ values to find how the nonlinear model is distinctive compared with the counterpart.  

\subsection{Strain $\epsilon_{22}$}

\begin{figure}[H]
	\centering
	\subfloat[\footnotesize{$\epsilon_{22}$ over the whole  surface.}]{\label{surffa}\includegraphics[width=0.5\textwidth]{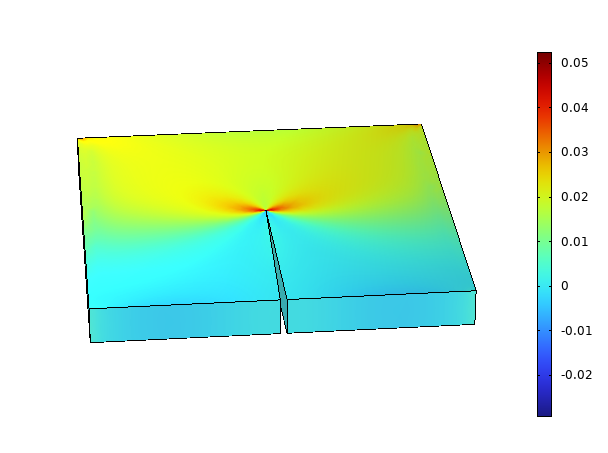}} 
	\subfloat[\footnotesize{$\epsilon_{22}$ over the slices parallel to the $yz$ plane. }]{\label{vslice}\includegraphics[width=0.5\textwidth]{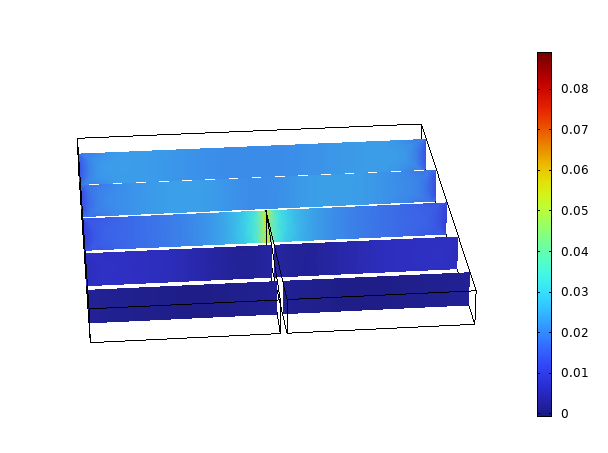}} \\
		\subfloat[\footnotesize{$\epsilon_{22}$ over the slices parallel to the $xz$ plane.}]{\label{xzslice}\includegraphics[width=0.5\textwidth]{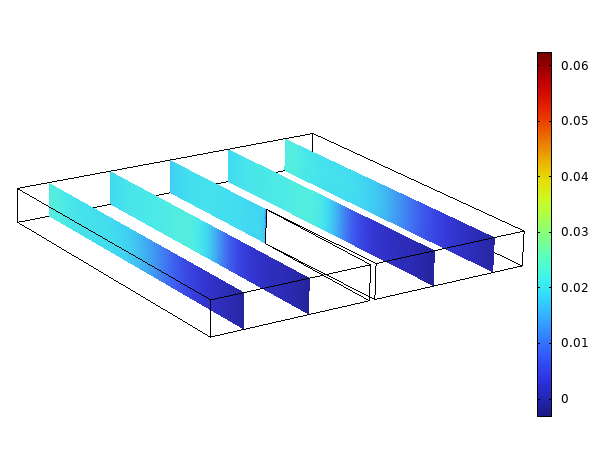}} 
	\subfloat[\footnotesize{$\epsilon_{22}$ over the  slice parallel to the $xy$ plane through the middle of the plate. }]{\label{xyslice}\includegraphics[width=0.5\textwidth]{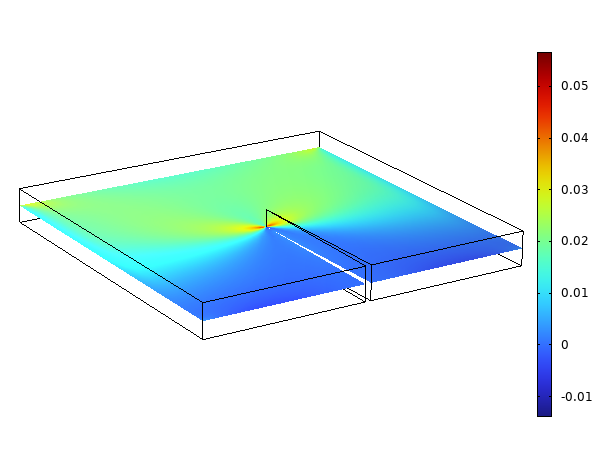}} 
	\caption{{\footnotesize{Illustration of $\epsilon_{22}$ over the whole geometrical surface and discrete slices in different directions. All the panels complement each other to facilitate understanding of the $\epsilon_{22}$ distribution. The concentration of $\epsilon_{22}$ forming two branches near the crack tip can be viewed particularly in Figs. \ref{surffa} and \ref{xyslice}.}}}
	\label{epsi22sli}
\end{figure}

 The general pattern for $\epsilon_{22}$ results is similar for all $\beta$ values used in the computation. We only illustrate the results for $\beta=-30$ as a representation. The 3-D results in the whole geometry, including the interior and surface, are impossible to display in one image. Instead, we show the distribution over the surface and discrete plane slices oriented in different directions (Fig. \ref{epsi22sli}). The $\epsilon_{22}$ distribution over the plate surface is displayed in Fig.  \ref{surffa}. In particular, the distribution over the top surface can be viewed.  The bottom surface displays similar results as the top surface but cannot be considered here. The results in the top surface show two branches of concentration of $\epsilon_{22}$ in red near the crack tip. Such concentration fades away towards the two back corners. Overall, over the top surface, the $\epsilon_{22}$ value is more significant in the back behind the crack tip than in the part in front of the crack tip. 

Figure \ref{vslice} shows the $\epsilon_{22}$ distribution over the slices parallel to the $yz$ plane.  The concentration of $\epsilon_{22}$ can be viewed in bright blue in the middle of the slice crossing the crack tip. Over each slice, the result shows an apparent left-right symmetry. If the slice is in front of the crack tip, $\epsilon_{22}$ is more diminutive, given by a darker color. On the contrary, $\epsilon_{22}$ is greater in slices in the back of the crack tip provided by a brighter color. Such results correspond to the outcomes in the top surface of the plate in Fig. \ref{surffa}. Figure \ref{xzslice} shows the $\epsilon_{22}$ distribution in the slices parallel to the $xz$ plane. Each slice except the middle one shows two distinct parts: dark and bright, transitioning through the crack tip. Lastly, Figure \ref{xyslice} shows the $\epsilon_{22}$ distribution over a single slice parallel to the $xy$ plane in the middle of the plate. The outcome is similar to the one on the top surface of the plate in Fig.  \ref{surffa} but bears different values. In summary, the four panels in Fig. \ref{epsi22sli} assist in viewing the 3-D $\epsilon_{22}$ distribution from different perspectives. Each panel's color scale shows a transparent color gradient for easy viewing of the $\epsilon_{22}$ distribution. 

\begin{figure}[H]
	\centering
	\subfloat[\footnotesize{For negative $\beta$ values.}]{\label{negabsmid}\includegraphics[width=0.5\textwidth]{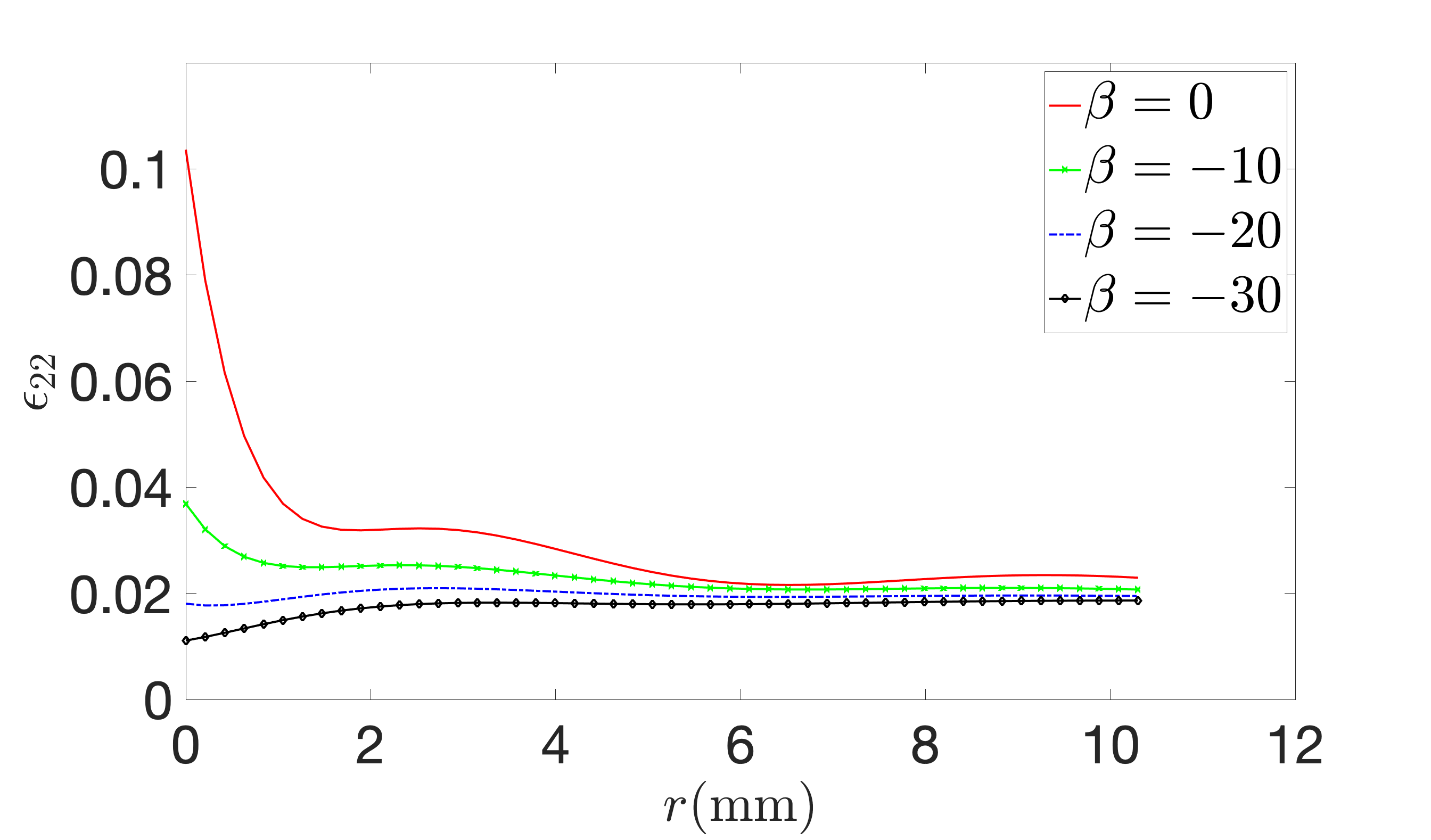}} 
	\subfloat[\footnotesize{For positive $\beta$ values.}]{\label{posbsmid}\includegraphics[width=0.5\textwidth]{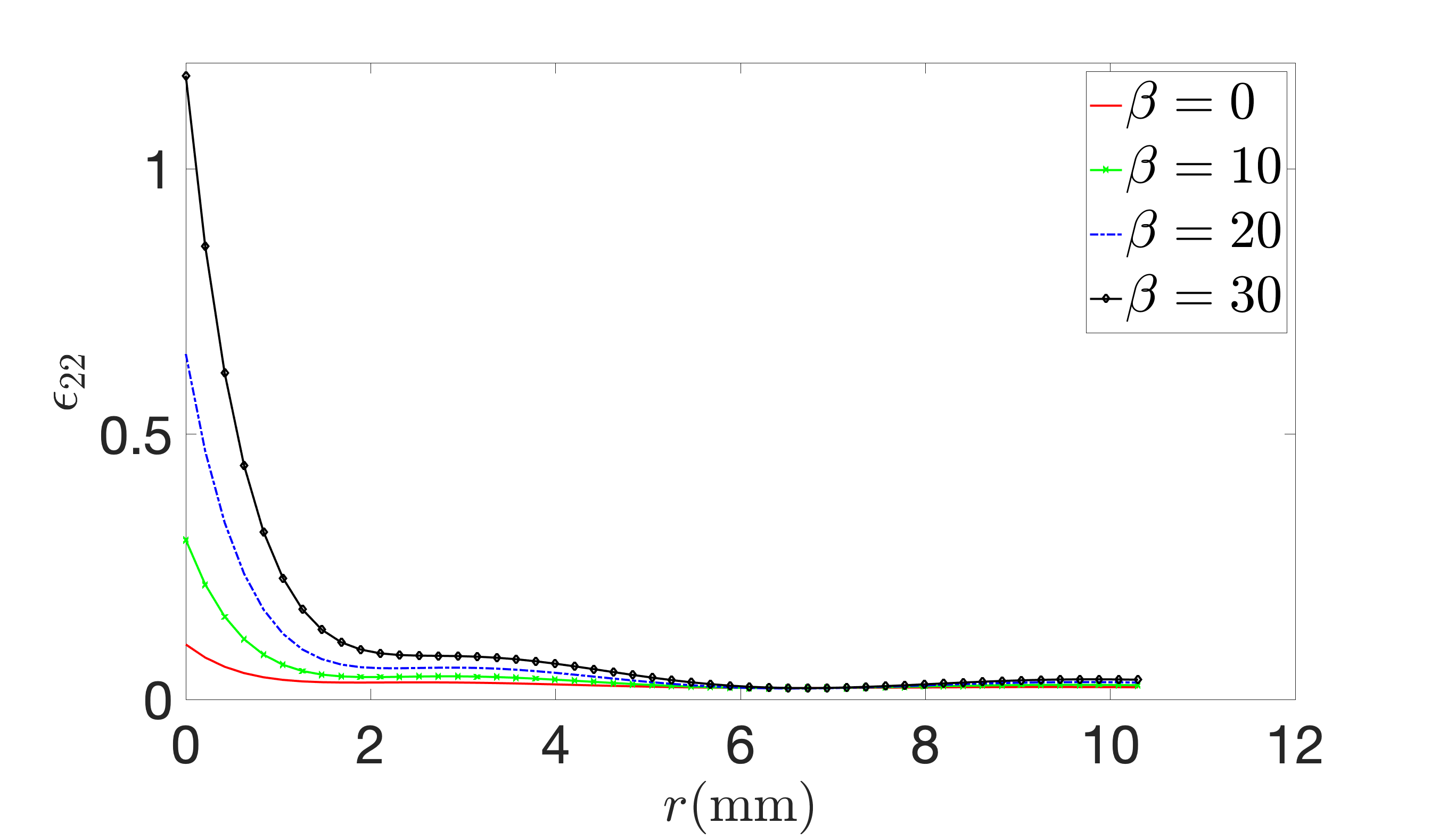}} 
	\caption{{\footnotesize{$\epsilon_{22}$ vs. $r$ for different $\beta$ values. For negative $\beta$, greater $|\beta|$ produces smaller $\epsilon_{22}$ at the crack-tip ($r=0$), realizing the strain-limiting effect. On the contrary, for positive $\beta$, greater $|\beta|$ produces greater $\epsilon_{22}$ at the crack-tip ($r=0$), failing to realize the strain-limit effect. All curves for both positive and negative $\beta$ values show slightly wavy behaviors.}}}
	\label{strainMidl}
\end{figure}

The relation between $\epsilon_{22}$ and $r$ is demonstrated in Fig. \ref{strainMidl}. The results for $\beta=0$ are for the linear-elasticity fracture model.  Smaller $\epsilon_{22}$ at the crack tip ($r=0$) than $\epsilon_{22}$ for $\beta=0$ is deemed to realize the strain-limiting effect. The negative and positive $\beta$ values results are shown in Fig. \ref{negabsmid} and Fig. \ref{posbsmid}, respectively.  Fig. \ref{negabsmid} shows that as $\beta$ is negative and decreases, $\epsilon_{22}$ at the crack-tip ($r=0$) decreases, reflecting the effect of the nonlinear model in reducing the strain near the crack tip. Fig. \ref{posbsmid} shows that as $\beta$ is positive and increases, greater $\beta$ generates greater $\epsilon_{22}$ at the crack-tip. The outcomes show that only negative $\beta$ values realize the nonlinear, strain-limiting effect and that positive $\beta$ values produce a counter-effect in this aspect. Furthermore, the slightly wavy behavior in the curves shows the impact of density-dependent moduli.

\subsection{Stress $T_{22}$}

\begin{figure}[H]
	\centering
	\includegraphics[width=0.5\textwidth]{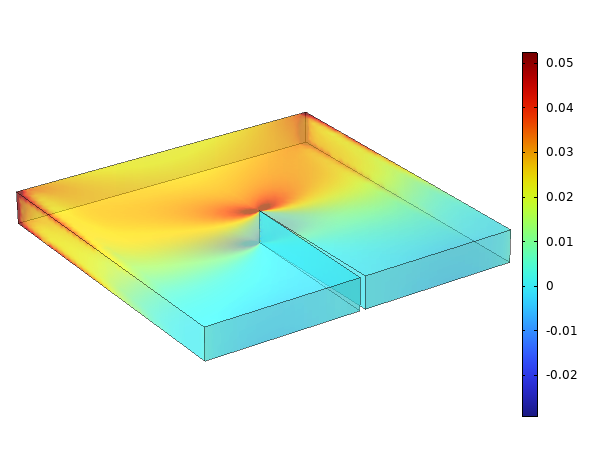}
	\caption{\footnotesize{A transparent surface view of $T_{22}$ distribution (unit: $10^4$ Pa) at $\beta=-30$. The concentration of $T_{22}$ can be viewed near the crack tip and on the rear part of the left and right lateral sides. The general pattern is similar to the one for $\epsilon_{22}$ in Fig. \ref{surffa}.} }
	\label{T22beta_30}
\end{figure}

The $T_{22}$ distribution over the plate is similar to the one for $\epsilon_{22}$ in Fig. \ref{surffa}.   We exhibit a transparent surface view of $T_{22}$ distribution at $\beta=-30$ in Fig. \ref{T22beta_30} as a representative for other $\beta$ values.  We skip showing the $T_{22}$ distribution over slices in different directions for brevity. In Fig. \ref{T22beta_30}, the concentration of $T_{22}$ can be viewed near the crack tip and on the rear half of the left and right lateral sides. Such concentration also spreads through a large area in the rear part of the plate. The distribution on the top and bottom surfaces is identical.  Such $T_{22}$ concentration is more intense than $\epsilon_{22}$ concentration, suggesting the stress is more sensitive than the strain under the impact of the parameter $\beta$.

\begin{figure}[H]
	\centering
	\subfloat[\footnotesize{For negative $\beta$ values.}]{\label{T22NegBeta}\includegraphics[width=0.5\textwidth]{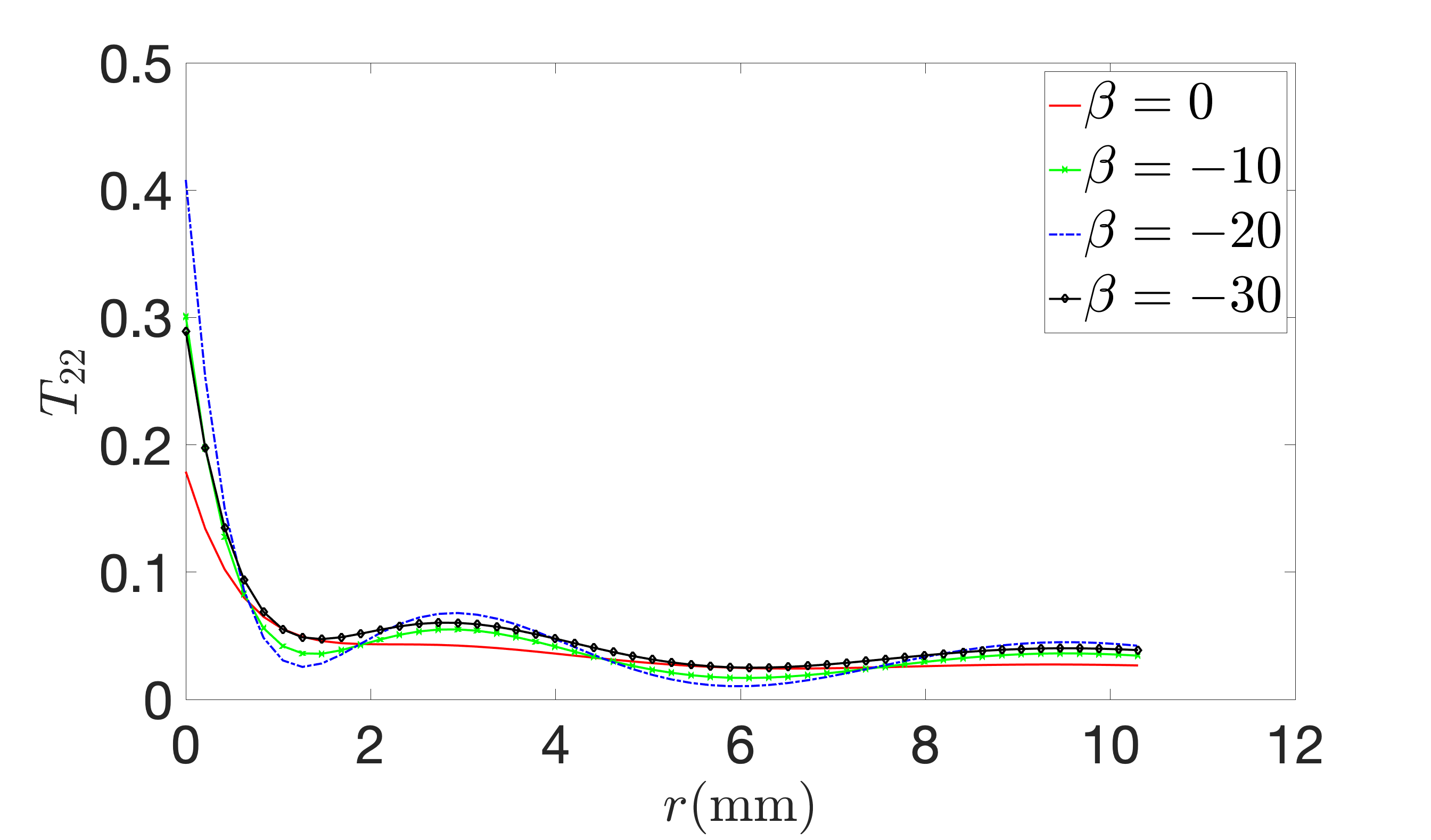}} 
	\subfloat[\footnotesize{For positive $\beta$ values.}]{\label{T22PosBeta}\includegraphics[width=0.5\textwidth]{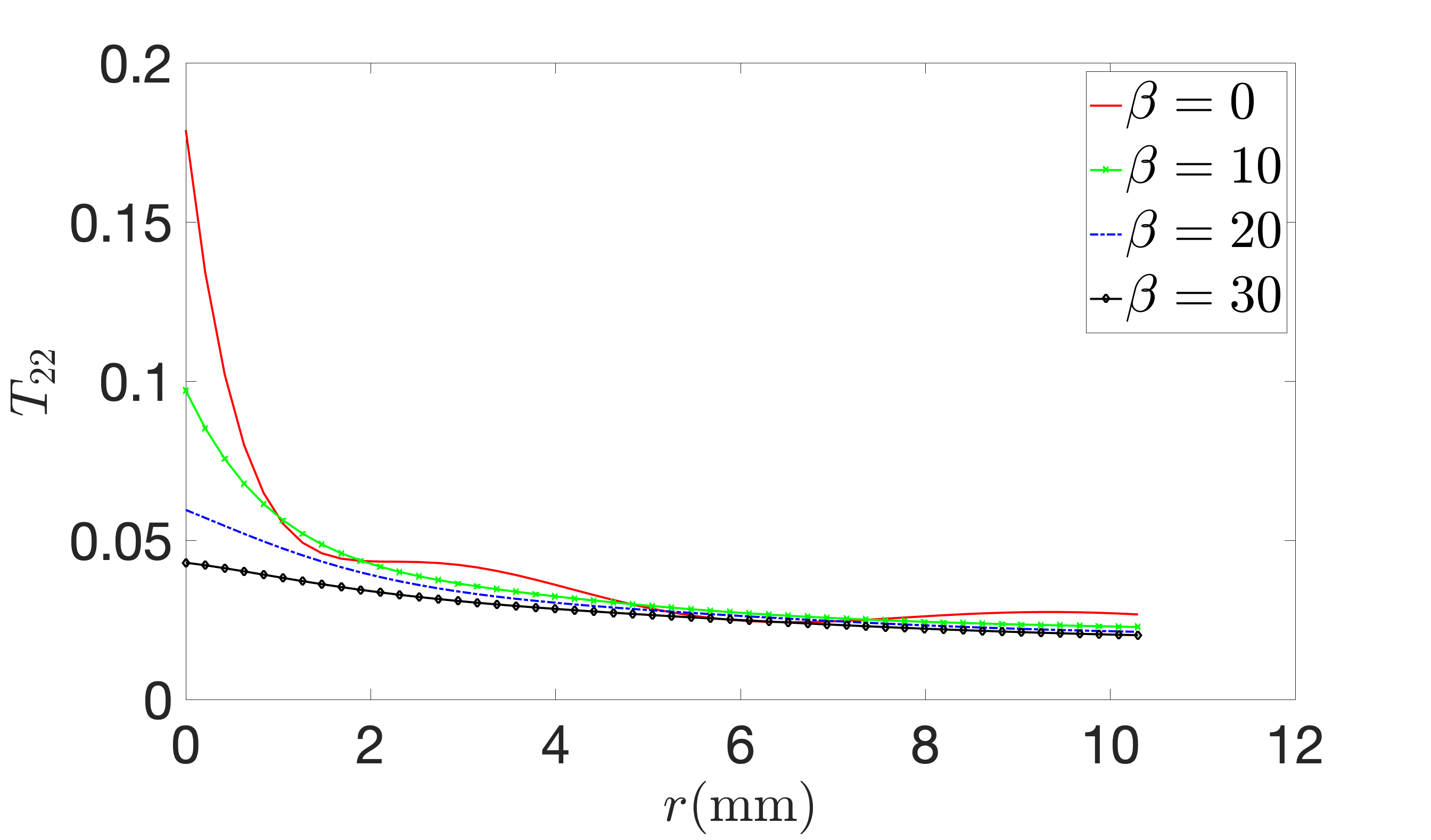}} 
	\caption{{\footnotesize{$T_{22}$ vs. $r$ ($T_{22}$ unit: $10^4$ Pa). In the first panel, at $r=0$ (the crack-tip), negative $\beta$ does not reduce $T_{22}$ compared with $T_{22}$ at $\beta=0$. Generally, a more negative value does not necessarily generate greater $T_{22}$ at the crack tip. In the second panel for $\beta>0$, positive $\beta$ generates smaller $T_{22}$ value compared with $T_{22}$ at $\beta=0$. Greater $\beta$ creates smaller $T_{22}$  at $r=0$. Slightly wavy patterns also appear in all the curves due to the effect of the density-dependent moduli.}}}
	\label{T22midl}
\end{figure}

The function of $T_{22}$ vs. $r$ is displayed in Fig. \ref{T22midl}. In Fig. \ref{T22NegBeta}, it shows negative $\beta$  fails to reduce $T_{22}$ at the crack tip for the linear-elasticity fracture model ($\beta=0$). On the contrary, negative $\beta$  increases $T_{22}$ at the crack tip, but a smaller negative $\beta$ value may not necessarily produce a larger increase for crack-tip $T_{22}$. This phenomenon is contrary to the $\beta$ effect on $\epsilon_{22}$ in Fig. \ref{negabsmid} for negative $\beta$. However, as shown in Fig. \ref{T22PosBeta},  positive $\beta$  is capable of reducing $T_{22}$ at the crack tip for $\beta=0$. The greater the positive $\beta$ is, the higher the reduction is. Similar to Fig. \ref{strainMidl}, all curves show slightly wavy patterns due to the impact of the density-dependent moduli.

\subsection{SIF}

{For the current problem, we do not have an asymptotic or analytical solution. Therefore, an explicit description of crack-tip SIF is not available. However, one can utilize the crack-tip SIF defined for the linear elastic fracture mechanics model and use the finite element solution for the nonlinear model to glean some vital physical insight into SIF.} \\

\begin{figure}[H]
	\centering
	\subfloat[\footnotesize{For negative $\beta$ values.}]{\label{negBetSIF}\includegraphics[width=0.5\textwidth]{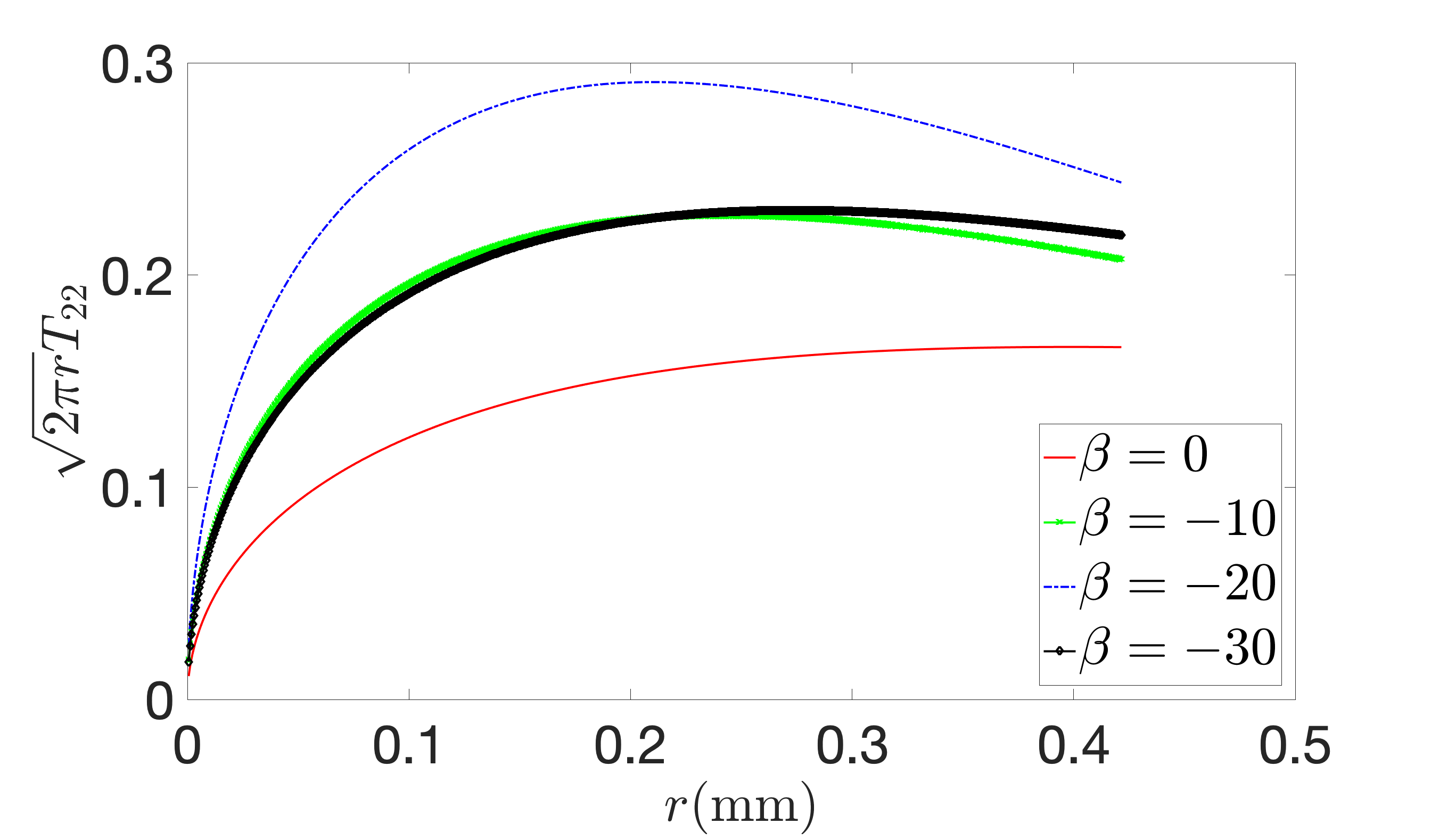}} 
	\subfloat[\footnotesize{For positive $\beta$ values.}]{\label{PosBetSIF}\includegraphics[width=0.5\textwidth]{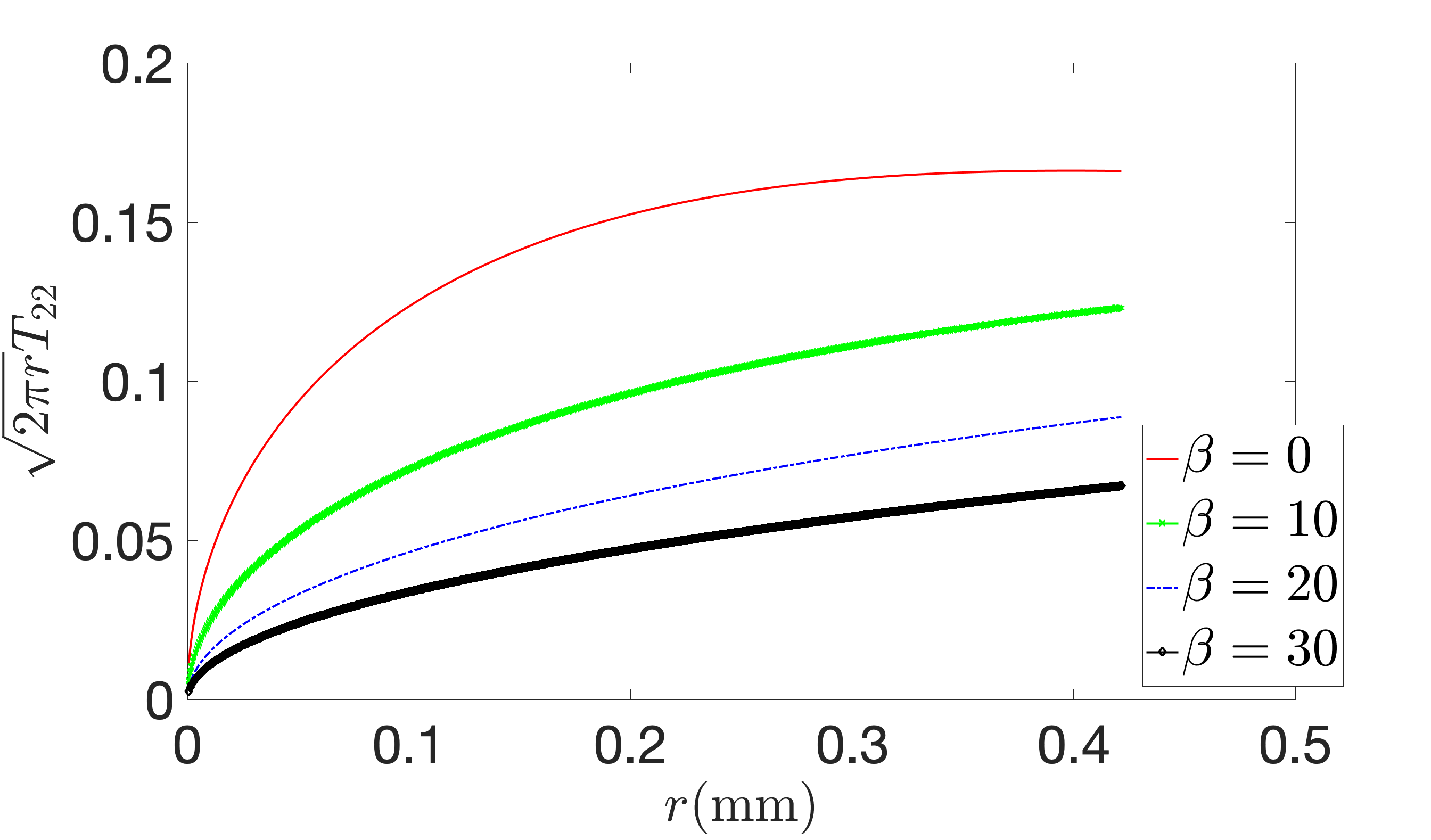}} 
	\caption{{\footnotesize{Exhibition of $\sqrt{2 \pi r} T_{22}$ vs. $r$ for studying the SIF (unit of $\sqrt{2 \pi r} T_{22}$: $10^4\text{mm}^{1/2}\text{Pa}$). All the curves approximately approach the same limit as $r$ approaches 0, implying a common SIF for all cases.}}}
	\label{SIFcur}
\end{figure}

The SIF $K_I$ in the context of linear elasticity  defined as
\begin{equation}\label{SIF_M1}
	K_I = \lim_{r \to 0^+} \sqrt{2 \pi r} T_{22}.
\end{equation}
It is challenging to analytically compute $K_I$ for complicated models such as the ones in 3-D that we are investigating in this paper. We compute $K_I$ numerically for illustrative purposes to understand its value as a function of the nonlinear modeling parameter $\beta$. The function $\sqrt{2 \pi r} T_{22}$ vs. $r$ is shown in Fig. \ref{SIFcur}. For $\beta=0, \,-10, \,-20, \,-30$, the $K_I$ values are, respectively, 0.0110, 0.0185, 0.0251, 0.0178; For $\beta=0, \,10, \,20, \,30$, the $K_I$ values are,  respectively, 0.0110, 0.0060, 0.0037, 0.0026 (Unit: $10^4\text{mm}^{1/2}\text{Pa}$). All the $K_I$ values are very close, implying an identical theoretical $K_I$ value for all $\beta$  is highly likely. Such numerical approximation for $K_I$ shows that the nonlinear model with density-dependent moduli is designed appropriately under the common crack criteria used for the linear-elasticity fracture model. Also, in the two panels of Fig. \ref{SIFcur}, the curve shapes for $\beta<0$ and $\beta>0$ are slightly distinct. The $\beta<0$ curves show a greater curvature, while those for $\beta>0$ show flatter behaviors. Such difference is caused by the sign of $\beta$, disclosing how negative and positive $\beta$ values impact the SIF differently.

It is clear from the above Fig~\ref{SIFcur} that the SIF values near the crack-tip are higher for the negative $\beta$ values; hence, there could be stress concentration only for negative $\beta$ values. 

\subsection{Strain energy}

\begin{figure}[H]
	\centering
	\subfloat[\footnotesize{For negative $\beta$ values.}]{\label{eneNegBe}\includegraphics[width=0.5\textwidth]{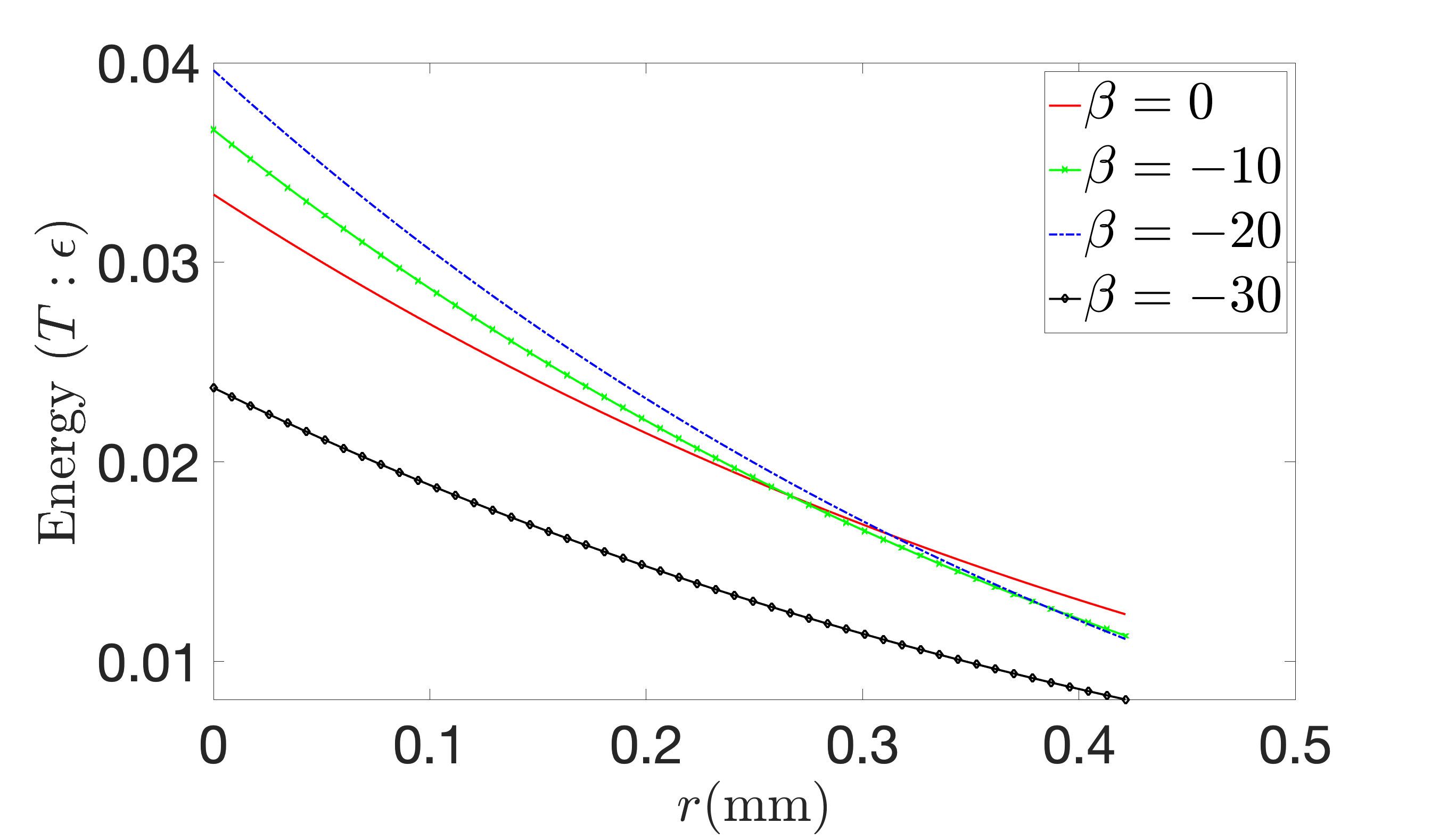}} 
	\subfloat[\footnotesize{For positive $\beta$ values.}]{\label{enePosBe}\includegraphics[width=0.5\textwidth]{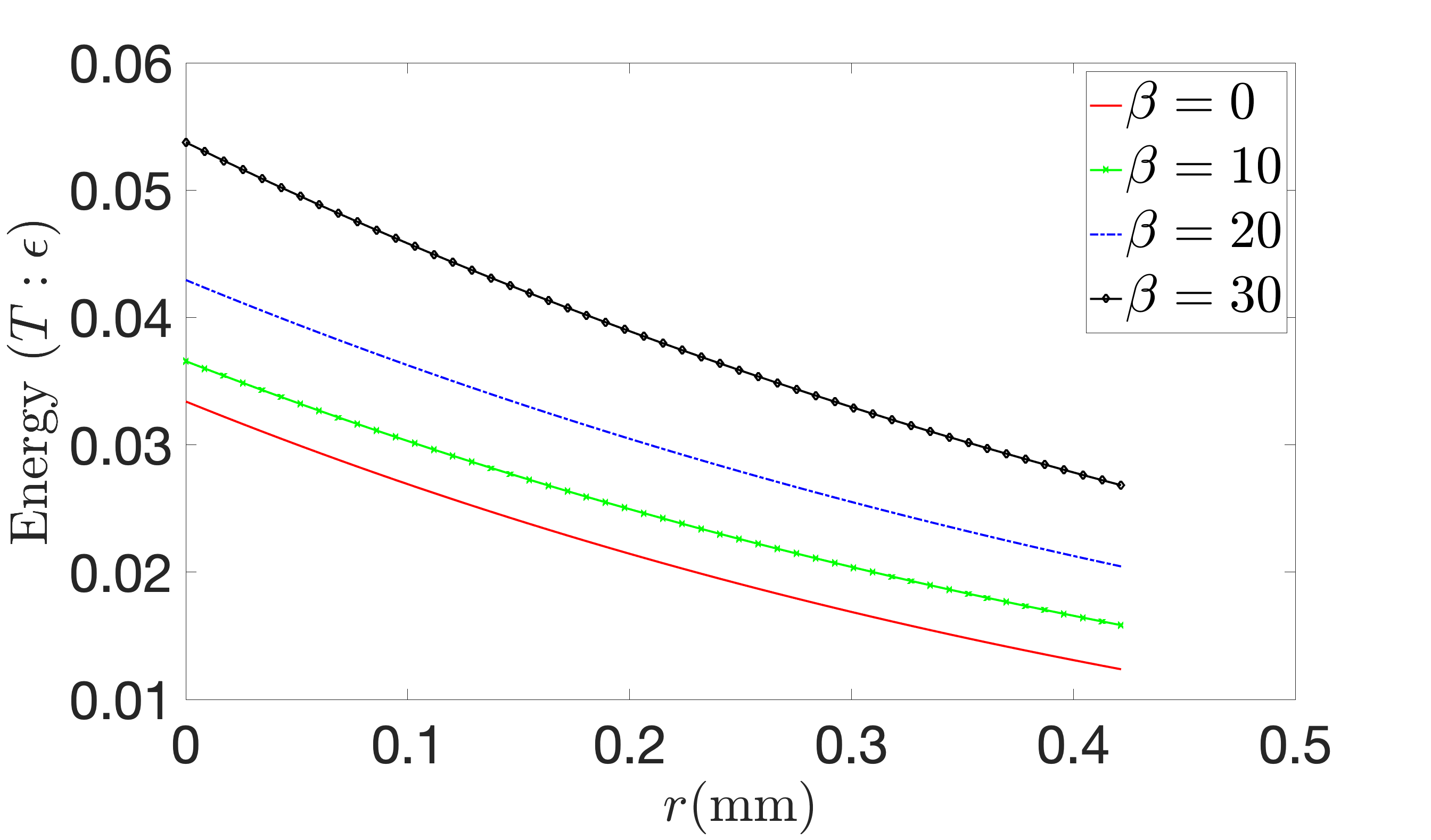}} 
	\caption{{\footnotesize{Energy vs. $r$ for both positive and negative $\beta$ values (unit of energy: $10^4\text{mm}^{1/2}\text{Pa}$). In the left panel, $\beta=-30$ gives the energy at the crack tip ($r=0$) smaller than that for $\beta=0$, and $\beta=-10$ and $\beta=-20$ generate greater energy than that for $\beta=0$ at $r=0$. In the right panel, all positive $\beta$ values generate energy at $r=0$ greater than that for $\beta=0$, and the greater $\beta$ is, the greater the energy at the crack tip is.}}}
	\label{energytip}
\end{figure}

The strain energy is defined to be $$(\mathbf{\epsilon}:\mathbf{T})=\sum_{i,j=1,2,3}\epsilon_{ij}T_{ij}.$$ We plot the function of strain energy vs. $r$ in Fig. \ref{energytip}. In Fig. \ref{eneNegBe}, $\beta=-30$ gives the smallest strain energy and particularly generates energy at the crack tip ($r=0$) smaller than that for $\beta=0$. However, both $\beta=-10$ and $\beta=-20$ generate greater energy than that for $\beta=0$ at $r=0$. Such outcomes suggest that not all negative $\beta$ values create smaller energy than that for the linear-elasticity fracture model.  Fig. \ref{enePosBe} shows that all positive $\beta$ values generate energy at $r=0$ greater than that for $\beta=0$. The pattern for energy under positive $\beta$ is that the greater the $\beta$ value is, the greater the energy at the crack tip is.

An important observation about the strain energy is that it gets its most significant value directly ahead of the crack tip, which is consistent with the observation made in the linearized elasticity model. Therefore, one can utilize the classical fracture criterion to study the quasi-static and dynamic propagation of crack tips within the framework of the models proposed in this paper. One can also employ the crack-regularization techniques such as the one proposed in \cite{yoon2021quasi,lee2022finite} to study the evolution problems in the context of nonlinear models.

\section{Conclusions and discussions}
In the last few years, there have been many theoretical investigations on developing algebraically nonlinear models for the description of geometrically linear as well as nonlinear elastic solids \cite{rajagopal2003,rajagopal2007elasticity,mai2015monotonicity,mai2015strong}. The ultimate goal of such studies is to formulate well-posed mathematical models for the bulk material behavior,  which can predict bounded crack-tip strain fields and preclude the necessity of using \textit{ad hoc} modeling augmentations such as cohesive or process zones to correct the erroneous predictions from linear elasticity. Perhaps the most crucial applicability of the modeling paradigm introduced in \cite{rajagopal2011modeling,gou2015modeling,Mallikarjunaiah2015,MalliPhD2015} is that it can pave a way to develop new fracture theories applicable to a wide variety of brittle materials.  This paper is the first part of the extended research program to develop theoretical and computational models for cracks and fractures in a new class of models describing the response of 3-D porous elastic solids. Using the particular subclass of constitutive relationship, we derive a quasilinear elliptic BVP for the static v-notch domain.

The 3-D problem is discretized using a continuous Galerkin-type finite element method coupled with the Picard-type linearization, meaning that the problem is linear at each iteration.  Hence, a unique solution is guaranteed, and a direct solver can be used to find the nodal values of the computational results. A few remarks about the numerical solution and its interpretation are as follows.
\begin{itemize}
\item[1.] We observe that the finite element discretization of the 3-D BVP is convergent, and Picard's iterations take a reasonable number of iterations to meet the stopping criterion. The relative tolerance 0.001 is employed as the stopping criterion with a maximum number of iterations of 1000.The initial guess is obtained by solving the BVP from the linearized elasticity, i.e., by taking $\beta=0$ in the nonlinear model. We notice that the numerical solution to the proposed nonlinear BVP converges to the corresponding one from the classical linear elasticity model as $\beta \to 0$.  
\item[2. ] For the tensile loading, with negative $\beta$ values, the nonlinear model predicts that the crack-tip strain growth is much smaller than the classical linear model. Our model displays a strain-limiting effect for negative $\beta$ values. However, a clear opposite pattern for the crack-tip strain values is observed for the positive $\beta$ values.  Such results are expected because the stiffness parameter $(1 + \beta \, \tr \, \bfeps)$ changes in the nonlinear model, different from the constant stiffness values in the classical linearized model. 
\item[3. ] The stress concentration increases if the negative $\beta$ value decreases; hence, the proposed model has the same effect on stress concentration as in the classical model. This implies that the proposed material model exhibits a stronger brittle effect for smaller negative  $\beta$ values. The largest value of $\bfT_{22}$ appears directly before the crack-tip, another important feature in the classical linear model. However, the stress concentration decreases in the nonlinear model with the increasing positive $\beta$ values. 
\item[4. ] As pointed out earlier, we currently do not have an asymptotic crack-tip solution for the proposed model; hence, the classical definition of SIF from linear elasticity is utilized to glean some vital information about the singularities. The finite element solution is used in the SIF formula to characterize the nature of SIF for the current model. The term $\sqrt{2\pi r} \bfT_{22}$ for the SIF approximately reaches a constant ahead of the crack-tip, and the values are higher for negative $\beta$ values, which indicates the possibility of stress singularity. For the model studied in this paper, the important open problem is to characterize the singularity via asymptotic analysis similar to \cite{gou2015modeling,rajagopal2011modeling}. 
\item[5. ]  The strain energy density also obtains its maximum before the crack tip. These results are consistent with those from the classical linear elasticity model. Hence, one can use the same local fracture criterion near the crack tip to study evolution problems within the framework of the nonlinear theory presented in this paper. The phase-field regularization, similar to the one done for strain-limiting models in \cite{yoon2021quasi,lee2022finite}, can also be used for the current framework to study quasi-static crack propagation. 
\item[6. ] The present study may incorporate surface mechanics as described in \cite{ferguson2015numerical}, thereby enabling the elimination of unphysical stress singularities. Consequently, a new local criterion predicated on cleavage stress or the maximum stress at the crack tip may be employed to examine crack-tip evolution issues. Another significant future direction is to examine the error analysis of the local discontinuous Galerkin-type finite element discretization \cite{manohar2024hp}. 
\end{itemize}

%
%
%
%
%
%
%

\bibliographystyle{unsrtnat}
\bibliography{references}

\begin{thebibliography}{62}
\providecommand{\natexlab}[1]{#1}
\providecommand{\url}[1]{\texttt{#1}}
\expandafter\ifx\csname urlstyle\endcsname\relax
  \providecommand{\doi}[1]{doi: #1}\else
  \providecommand{\doi}{doi: \begingroup \urlstyle{rm}\Url}\fi

\bibitem[Rubinstein and Panyukov(2002)]{rubinstein2002elasticity}
M.~Rubinstein and S.~Panyukov.
\newblock Elasticity of polymer networks.
\newblock \emph{Macromolecules}, 35\penalty0 (17):\penalty0 6670--6686, 2002.

\bibitem[Leventis et~al.(2002)Leventis, Sotiriou-Leventis, Zhang, and
  Rawashdeh]{leventis2002nanoengineering}
N.~Leventis, C.~Sotiriou-Leventis, G.~Zhang, and A.-M.~M. Rawashdeh.
\newblock Nanoengineering strong silica aerogels.
\newblock \emph{Nano letters}, 2\penalty0 (9):\penalty0 957--960, 2002.

\bibitem[Chandrasekaran et~al.(2017)Chandrasekaran, Campbell, Baumann, and
  Worsley]{chandrasekaran2017carbon}
S.~Chandrasekaran, P.~G. Campbell, T.~F. Baumann, and M.~A. Worsley.
\newblock Carbon aerogel evolution: Allotrope, graphene-inspired, and
  3d-printed aerogels.
\newblock \emph{Journal of Materials Research}, 32\penalty0 (22):\penalty0
  4166--4185, 2017.

\bibitem[Anjos et~al.(2010)Anjos, Pereira, and Rosa]{anjos2010tensile}
O.~Anjos, H.~Pereira, and M.~E. Rosa.
\newblock Tensile properties of cork in axial stress and influence of porosity,
  density, quality and radial position in the plank.
\newblock \emph{European Journal of Wood and Wood Products}, 69\penalty0
  (1):\penalty0 85--91, 2010.

\bibitem[Alagappan et~al.(2023)Alagappan, Arumugam, and
  Rajagopal]{ALAGAPPAN2023100162}
P.~Alagappan, J.~Arumugam, and K.R. Rajagopal.
\newblock A note on the response of elastic bodies whose material moduli depend
  on the density and the mechanical pressure.
\newblock \emph{Applications in Engineering Science}, 16:\penalty0 100162,
  2023.
\newblock ISSN 2666-4968.
\newblock \doi{https://doi.org/10.1016/j.apples.2023.100162}.

\bibitem[Rajagopal(2003)]{rajagopal2003}
K~R Rajagopal.
\newblock On implicit constitutive theories.
\newblock \emph{Applications of Mathematics}, 48\penalty0 (4):\penalty0
  279--319, 2003.

\bibitem[Rajagopal(2007)]{rajagopal2007elasticity}
K.~R. Rajagopal.
\newblock The elasticity of elasticity.
\newblock \emph{Zeitschrift f{\"u}r Angewandte Mathematik und Physik (ZAMP)},
  58\penalty0 (2):\penalty0 309--317, 2007.

\bibitem[Rajagopal(2014)]{rajagopal2014nonlinear}
K.~R. Rajagopal.
\newblock On the nonlinear elastic response of bodies in the small strain
  range.
\newblock \emph{Acta Mechanica}, 225\penalty0 (6):\penalty0 1545--1553, 2014.

\bibitem[Rajagopal and Srinivasa(2007)]{rajagopal2007response}
K.~R. Rajagopal and A.~R. Srinivasa.
\newblock On the response of non-dissipative solids.
\newblock \emph{Proceedings of the Royal Society A: Mathematical, Physical and
  Engineering Sciences}, 463\penalty0 (2078):\penalty0 357--367, 2007.

\bibitem[Rajagopal and Srinivasa(2009)]{rajagopal2009class}
K.~R. Rajagopal and A.~R. Srinivasa.
\newblock On a class of non-dissipative materials that are not hyperelastic.
\newblock \emph{Proceedings of the Royal Society A: Mathematical, Physical and
  Engineering Sciences}, 465\penalty0 (2102):\penalty0 493--500, 2009.

\bibitem[Bridges and Rajagopal(2015)]{bridges2015implicit}
C.~Bridges and K.~R. Rajagopal.
\newblock Implicit constitutive models with a thermodynamic basis: a study of
  stress concentration.
\newblock \emph{Zeitschrift f{\"u}r angewandte Mathematik und Physik},
  66:\penalty0 191--208, 2015.

\bibitem[Mallikarjunaiah and Walton(2015)]{Mallikarjunaiah2015}
S.~M. Mallikarjunaiah and J.~R. Walton.
\newblock On the direct numerical simulation of plane-strain fracture in a
  class of strain-limiting anisotropic elastic bodies.
\newblock \emph{International Journal of Fracture}, 192\penalty0 (2):\penalty0
  217--232, Apr 2015.

\bibitem[Mallikarjunaiah(2015)]{MalliPhD2015}
S.~M. Mallikarjunaiah.
\newblock \emph{On Two Theories for Brittle Fracture: Modeling and Direct
  Numerical Simulations}.
\newblock PhD thesis, Texas A\&M University, 2015.

\bibitem[Ortiz-Bernardin et~al.(2014)Ortiz-Bernardin, Bustamante, and
  Rajagopal]{ortiz2014numerical}
A.~Ortiz-Bernardin, R.~Bustamante, and K.~R. Rajagopal.
\newblock A numerical study of elastic bodies that are described by
  constitutive equations that exhibit limited strains.
\newblock \emph{International Journal of Solids and Structures}, 51\penalty0
  (3-4):\penalty0 875--885, 2014.

\bibitem[Ortiz-Bernardin et~al.(2012)Ortiz-Bernardin, Bustamante, and
  Rajagopal]{ortiz2012}
A.~Ortiz-Bernardin, R.~Bustamante, and K.~R. Rajagopal.
\newblock A numerical study of a plate with a hole for a new class of elastic
  bodies.
\newblock \emph{Acta Mechanica}, 223\penalty0 (9):\penalty0 1971--1981, 2012.

\bibitem[Muliana et~al.(2018)Muliana, Rajagopal, Tscharnuter, Schrittesser, and
  Saccomandi]{muliana2018determining}
A.~Muliana, K.~R. Rajagopal, D.~Tscharnuter, B.~Schrittesser, and
  G.~Saccomandi.
\newblock Determining material properties of natural rubber using fewer
  material moduli in virtue of a novel constitutive approach for elastic
  bodies.
\newblock \emph{Rubber Chemistry and Technology}, 91\penalty0 (2):\penalty0
  375--389, 2018.

\bibitem[Kowalczyk-Gajewska et~al.(2019)Kowalczyk-Gajewska, Pieczyska,
  Golasi{\'n}ski, Maj, Kuramoto, and Furuta]{kowalczyk2019finite}
K.~Kowalczyk-Gajewska, E.~A. Pieczyska, K.~Golasi{\'n}ski, M.~Maj, S.~Kuramoto,
  and T.~Furuta.
\newblock A finite strain elastic-viscoplastic model of gum metal.
\newblock \emph{International Journal of Plasticity}, 119:\penalty0 85--101,
  2019.

\bibitem[Bustamante et~al.(2020)Bustamante, Montero, and
  Ortiz-Bernardin]{bustamante2020novel}
R.~Bustamante, S.~Montero, and A.~Ortiz-Bernardin.
\newblock A novel nonlinear constitutive model for rock: Numerical assessment
  and benchmarking.
\newblock \emph{Applications in Engineering Science}, 3:\penalty0 100012, 2020.

\bibitem[Bustamante and Rajagopal(2021)]{bustamante2021new}
R.~Bustamante and K.~R. Rajagopal.
\newblock A new type of constitutive equation for nonlinear elastic bodies.
  fitting with experimental data for rubber-like materials.
\newblock \emph{Proceedings of the Royal Society A}, 477\penalty0
  (2252):\penalty0 20210330, 2021.

\bibitem[Ogden(1972)]{ogden1972large}
R.~W. Ogden.
\newblock Large deformation isotropic elasticity--on the correlation of theory
  and experiment for incompressible rubberlike solids.
\newblock \emph{Proceedings of the Royal Society of London. A. Mathematical and
  Physical Sciences}, 326\penalty0 (1567):\penalty0 565--584, 1972.

\bibitem[Saito et~al.(2003)Saito, Furuta, Hwang, Kuramoto, Nishino, Suzuki,
  Chen, Yamada, Ito, Seno, Nonaka, Ikehata, Nagasako, Iwamoto, Ikuhara, and
  Sakuma]{saito2003multifunctional}
T.~Saito, T.~Furuta, J.-H. Hwang, S.~Kuramoto, K.~Nishino, N.~Suzuki, R.~Chen,
  A.~Yamada, K.~Ito, Y.~Seno, T.~Nonaka, H.~Ikehata, N.~Nagasako, C.~Iwamoto,
  Y.~Ikuhara, and T.~Sakuma.
\newblock Multifunctional alloys obtained via a dislocation-free plastic
  deformation mechanism.
\newblock \emph{Science}, 300\penalty0 (5618):\penalty0 464--467, 2003.
\newblock ISSN 0036-8075.
\newblock \doi{10.1126/science.1081957}.

\bibitem[Withey et~al.(2008)Withey, Jin, Minor, Kuramoto, Chrzan, and
  Morris~Jr]{withey2008deformation}
E.~Withey, M.~Jin, A.~Minor, S.~Kuramoto, D.~C. Chrzan, and J.~W. Morris~Jr.
\newblock The deformation of “gum metal” in nanoindentation.
\newblock \emph{Materials Science and Engineering: A}, 493\penalty0
  (1-2):\penalty0 26--32, 2008.

\bibitem[Hou et~al.(2009)Hou, Li, Hao, and Yang]{zhang2009fatigue}
F.~Q. Hou, S.~J. Li, Y.~L. Hao, and R~Yang.
\newblock Fatigue properties of a multifunctional titanium alloy exhibiting
  nonlinear elastic deformation behavior.
\newblock \emph{Scripta Materialia}, 60\penalty0 (8):\penalty0 733--736, 2009.

\bibitem[Johnson and Rasolofosaon(1996)]{johnson1996manifestation}
P.~A. Johnson and P.~N.~J. Rasolofosaon.
\newblock Manifestation of nonlinear elasticity in rock: convincing evidence
  over large frequency and strain intervals from laboratory studies.
\newblock \emph{Nonlinear processes in geophysics}, 3\penalty0 (2):\penalty0
  77--88, 1996.

\bibitem[Kulvait et~al.(2019)Kulvait, M{\'a}lek, and
  Rajagopal]{kulvait2019state}
V.~Kulvait, J.~M{\'a}lek, and K.~R. Rajagopal.
\newblock The state of stress and strain adjacent to notches in a new class of
  nonlinear elastic bodies.
\newblock \emph{Journal of Elasticity}, 135\penalty0 (1):\penalty0 375--397,
  2019.

\bibitem[Tian et~al.(2015)Tian, Yu, Ong, and Cui]{tian2015nonlinear}
Y.~Tian, Z.~Yu, C.~Y.~A. Ong, and W.~Cui.
\newblock Nonlinear elastic behavior induced by nano-scale $\alpha$ phase in
  $\beta$ matrix of $\beta$-type {T}i--25{N}b--3{Z}r--2{S}n--3{M}o titanium
  alloy.
\newblock \emph{Materials Letters}, 145:\penalty0 283--286, 2015.

\bibitem[Devendiran et~al.(2017)Devendiran, Sandeep, Kannan, and
  Rajagopal]{devendiran2017thermodynamically}
V.~K. Devendiran, R.~K. Sandeep, K.~Kannan, and K.~R. Rajagopal.
\newblock A thermodynamically consistent constitutive equation for describing
  the response exhibited by several alloys and the study of a meaningful
  physical problem.
\newblock \emph{International Journal of Solids and Structures}, 108:\penalty0
  1--10, 2017.

\bibitem[Hao et~al.(2005)Hao, Li, Sun, Zheng, Hu, and Yang]{hao2005super}
Y.~L. Hao, S.~J. Li, S.~Y. Sun, C.~Y. Zheng, Q.~M. Hu, and R.~Yang.
\newblock Super-elastic titanium alloy with unstable plastic deformation.
\newblock \emph{Applied Physics Letters}, 87\penalty0 (9):\penalty0 091906,
  2005.

\bibitem[Rajagopal and Walton(2011)]{rajagopal2011modeling}
K.~R. Rajagopal and J.~R. Walton.
\newblock Modeling fracture in the context of a strain-limiting theory of
  elasticity: a single anti-plane shear crack.
\newblock \emph{International journal of fracture}, 169\penalty0 (1):\penalty0
  39--48, 2011.

\bibitem[Gou et~al.(2015)Gou, Mallikarjunaiah, Rajagopal, and
  Walton]{gou2015modeling}
K.~Gou, S.~M. Mallikarjunaiah, K.~R. Rajagopal, and J.~R. Walton.
\newblock Modeling fracture in the context of a strain-limiting theory of
  elasticity: A single plane-strain crack.
\newblock \emph{International Journal of Engineering Science}, 88:\penalty0
  73--82, 2015.

\bibitem[Yoon and Mallikarjunaiah(2022)]{HCY_SMM_MMS2022}
H.~C. Yoon and S.~M. Mallikarjunaiah.
\newblock A finite-element discretization of some boundary value problems for
  nonlinear strain-limiting elastic bodies.
\newblock \emph{Mathematics and Mechanics of Solids}, 27\penalty0 (2):\penalty0
  281--307, 2022.

\bibitem[Yoon et~al.(2021)Yoon, Lee, and Mallikarjunaiah]{yoon2021quasi}
H.~C. Yoon, S.~Lee, and S.~M. Mallikarjunaiah.
\newblock Quasi-static anti-plane shear crack propagation in nonlinear
  strain-limiting elastic solids using phase-field approach.
\newblock \emph{International Journal of Fracture}, 227\penalty0 (2):\penalty0
  153--172, 2021.

\bibitem[Lee et~al.(2022)Lee, Yoon, and Mallikarjunaiah]{lee2022finite}
S.~Lee, H.~C. Yoon, and S.~M. Mallikarjunaiah.
\newblock Finite element simulation of quasi-static tensile fracture in
  nonlinear strain-limiting solids with the phase-field approach.
\newblock \emph{Journal of Computational and Applied Mathematics},
  399:\penalty0 113715, 2022.

\bibitem[Yoon et~al.(2022)Yoon, Vasudeva, and Mallikarjunaiah]{yoon2022finite}
H.~C. Yoon, K.~K. Vasudeva, and S.~M. Mallikarjunaiah.
\newblock Finite element model for a coupled thermo-mechanical system in
  nonlinear strain-limiting thermoelastic body.
\newblock \emph{Communications in Nonlinear Science and Numerical Simulation},
  page 106262, 2022.

\bibitem[Gou and Mallikarjunaiah(2023{\natexlab{a}})]{gou2023MMS}
K.~Gou and S.~M. Mallikarjunaiah.
\newblock Finite element study of v-shaped crack-tip fields in a
  three-dimensional nonlinear strain-limiting elastic body.
\newblock \emph{Mathematicas and Mechanics of Solids}, pages 1--22,
  2023{\natexlab{a}}.

\bibitem[Gou and Mallikarjunaiah(2023{\natexlab{b}})]{gou2023computational}
K.~Gou and S.~M. Mallikarjunaiah.
\newblock Computational modeling of circular crack-tip fields under tensile
  loading in a strain-limiting elastic solid.
\newblock \emph{Communications in Nonlinear Science and Numerical Simulation},
  121:\penalty0 107217, 2023{\natexlab{b}}.

\bibitem[Rajagopal(2021{\natexlab{a}})]{rajagopal2021b}
K.~R. Rajagopal.
\newblock An implicit constitutive relation for describing the small strain
  response of porous elastic solids whose material moduli are dependent on the
  density.
\newblock \emph{Mathematics and Mechanics of Solids}, page 10812865211021465,
  2021{\natexlab{a}}.

\bibitem[Rajagopal and Saccomandi(2022)]{rajagopal2022implicit}
K.~R. Rajagopal and G.~Saccomandi.
\newblock Implicit nonlinear elastic bodies with density dependent material
  moduli and its linearization.
\newblock \emph{International Journal of Solids and Structures}, 234:\penalty0
  111255, 2022.

\bibitem[Itou et~al.(2022)Itou, Kovtunenko, and
  Rajagopal]{itou2022investigation}
H.~Itou, V.~A. Kovtunenko, and K.~R. Rajagopal.
\newblock Investigation of implicit constitutive relations in which both the
  stress and strain appear linearly, adjacent to non-penetrating cracks.
\newblock \emph{Mathematical Models and Methods in Applied Sciences},
  32\penalty0 (07):\penalty0 1475--1492, 2022.

\bibitem[Itou et~al.(2021)Itou, Kovtunenko, and Rajagopal]{itou2021implicit}
H.~Itou, V.~A. Kovtunenko, and K.~R. Rajagopal.
\newblock On an implicit model linear in both stress and strain to describe the
  response of porous solids.
\newblock \emph{Journal of Elasticity}, 144\penalty0 (1):\penalty0 107--118,
  2021.

\bibitem[Itou et~al.(2023)Itou, Kovtunenko, and
  Rajagopal]{itou2023generalization}
H.~Itou, V.~A. Kovtunenko, and K.~R. Rajagopal.
\newblock A generalization of the kelvin--voigt model with pressure-dependent
  moduli in which both stress and strain appear linearly.
\newblock \emph{Mathematical Methods in the Applied Sciences}, 2023.

\bibitem[Pr{\u{u}}{\v{s}}a et~al.(2022)Pr{\u{u}}{\v{s}}a, Rajagopal, and
  Wineman]{pruuvsa2022pure}
V.~Pr{\u{u}}{\v{s}}a, K.~R. Rajagopal, and A.~Wineman.
\newblock Pure bending of an elastic prismatic beam made of a material with
  density-dependent material parameters.
\newblock \emph{Mathematics and Mechanics of Solids}, 27\penalty0 (8):\penalty0
  1546--1558, 2022.

\bibitem[Pr{\u{u}}{\v{s}}a and Trnka(2023)]{pruuvsa2023mechanical}
V.~Pr{\u{u}}{\v{s}}a and L.~Trnka.
\newblock Mechanical response of elastic materials with density-dependent young
  modulus.
\newblock \emph{Applications in Engineering Science}, 14:\penalty0 100126,
  2023.

\bibitem[Yoon et~al.(2023)Yoon, Mallikarjunaiah, and Bhatta]{HYCSMM2023}
H.C. Yoon, S.~M. Mallikarjunaiah, and D.~Bhatta.
\newblock Finite element solution of crack-tip fields for an elastic porous
  solid with density-dependent material moduli and preferential stiffness.
\newblock \emph{Under review}, 2023.

\bibitem[Rajagopal(2021{\natexlab{b}})]{rajagopal2021a}
K.~R. Rajagopal.
\newblock An implicit constitutive relation in which the stress and the
  linearized strain appear linearly, for describing the small displacement
  gradient response of elastic solids.
\newblock \emph{arXiv preprint arXiv:2101.01208}, 2021{\natexlab{b}}.

\bibitem[Gokulnath et~al.(2017)Gokulnath, Saravanan, and
  Rajagopal]{gokulnath2017representations}
C.~Gokulnath, U.~Saravanan, and K.~R. Rajagopal.
\newblock Representations for implicit constitutive relations describing
  non-dissipative response of isotropic materials.
\newblock \emph{Zeitschrift f{\"u}r angewandte Mathematik und Physik},
  68:\penalty0 1--14, 2017.

\bibitem[Murru and Rajagopal(2021{\natexlab{a}})]{murru2021stress}
P.~T. Murru and K.~R. Rajagopal.
\newblock Stress concentration due to the presence of a hole within the context
  of elastic bodies.
\newblock \emph{Material Design \& Processing Communications}, page e219,
  2021{\natexlab{a}}.

\bibitem[Murru and Rajagopal(2021{\natexlab{b}})]{murru2021ZAMM}
P.~Murru and K.~R. Rajagopal.
\newblock Stress concentration due to the bi-axial deformation of a plate of a
  porous elastic body with a hole.
\newblock \emph{ZAMM-Journal of Applied Mathematics and Mechanics/Zeitschrift
  f{\"u}r Angewandte Mathematik und Mechanik}, 101\penalty0 (11):\penalty0
  e202100103, 2021{\natexlab{b}}.

\bibitem[Vajipeyajula et~al.(2023)Vajipeyajula, Murru, and
  Rajagopal]{vajipeyajula2023stress}
B.~Vajipeyajula, P.~Murru, and K.~R. Rajagopal.
\newblock Stress concentration due to an elliptic hole in a porous elastic
  plate.
\newblock \emph{Mathematics and Mechanics of Solids}, 28\penalty0 (3):\penalty0
  854--869, 2023.

\bibitem[Murru et~al.(2022)Murru, Torrence, Grasley, Rajagopal, Alagappan, and
  Garboczi]{murru2022density}
P.~T. Murru, C.~Torrence, Z.~Grasley, K.~R. Rajagopal, P.~Alagappan, and
  E.~Garboczi.
\newblock Density-driven damage mechanics (d3-m) model for concrete i:
  mechanical damage.
\newblock \emph{International Journal of Pavement Engineering}, 23\penalty0
  (4):\penalty0 1161--1174, 2022.

\bibitem[Knees et~al.(2015)Knees, Rossi, and Zanini]{knees2015quasilinear}
D.~Knees, R.~Rossi, and C.~Zanini.
\newblock A quasilinear differential inclusion for viscous and rate-independent
  damage systems in non-smooth domains.
\newblock \emph{Nonlinear Analysis: Real World Applications}, 24:\penalty0
  126--162, 2015.

\bibitem[Negri and Scala(2017)]{negri2017quasi}
M.~Negri and R.~Scala.
\newblock A quasi-static evolution generated by local energy minimizers for an
  elastic material with a cohesive interface.
\newblock \emph{Nonlinear Analysis: Real World Applications}, 38:\penalty0
  271--305, 2017.

\bibitem[Kuttler and Shillor(2006)]{kuttler2006quasistatic}
K.~L. Kuttler and M.~Shillor.
\newblock Quasistatic evolution of damage in an elastic body.
\newblock \emph{Nonlinear analysis: real world applications}, 7\penalty0
  (4):\penalty0 674--699, 2006.

\bibitem[Sendova and Walton(2010)]{sendova2010new}
T.~Sendova and J.~R. Walton.
\newblock A new approach to the modeling and analysis of fracture through
  extension of continuum mechanics to the nanoscale.
\newblock \emph{Mathematics and Mechanics of Solids}, 15\penalty0 (3):\penalty0
  368--413, 2010.

\bibitem[Walton(2012)]{walton2012note}
J.~R. Walton.
\newblock A note on fracture models incorporating surface elasticity.
\newblock \emph{Journal of Elasticity}, 109:\penalty0 95--102, 2012.

\bibitem[Ferguson et~al.(2015)Ferguson, Muddamallappa, and
  Walton]{ferguson2015numerical}
L.~A. Ferguson, M.~Muddamallappa, and J.~R. Walton.
\newblock Numerical simulation of mode-iii fracture incorporating interfacial
  mechanics.
\newblock \emph{International Journal of Fracture}, 192:\penalty0 47--56, 2015.

\bibitem[Sinclair(2004)]{sinclair2004stress}
G.~B. Sinclair.
\newblock Stress singularities in classical elasticity--i: Removal,
  interpretation, and analysis.
\newblock \emph{Appl. Mech. Rev.}, 57\penalty0 (4):\penalty0 251--298, 2004.

\bibitem[Kim et~al.(2011)Kim, Schiavone, and Ru]{kim2011analysis}
C.~I. Kim, P.~Schiavone, and C.-Q. Ru.
\newblock Analysis of plane-strain crack problems (mode-i \& mode-ii) in the
  presence of surface elasticity.
\newblock \emph{Journal of Elasticity}, 104:\penalty0 397--420, 2011.

\bibitem[Ciarlet(2002)]{ciarlet2002finite}
P.~G. Ciarlet.
\newblock \emph{The Finite Element Method for Elliptic Problems}.
\newblock SIAM, 2002.

\bibitem[Mai and Walton(2015{\natexlab{a}})]{mai2015monotonicity}
T.~Mai and J.~R. Walton.
\newblock On monotonicity for strain-limiting theories of elasticity.
\newblock \emph{Journal of Elasticity}, 120:\penalty0 39--65,
  2015{\natexlab{a}}.

\bibitem[Mai and Walton(2015{\natexlab{b}})]{mai2015strong}
T.~Mai and J.~R. Walton.
\newblock On strong ellipticity for implicit and strain-limiting theories of
  elasticity.
\newblock \emph{Mathematics and Mechanics of Solids}, 20\penalty0 (2):\penalty0
  121--139, 2015{\natexlab{b}}.

\bibitem[Manohar and Mallikarjunaiah(2024)]{manohar2024hp}
R.~Manohar and S.~M. Mallikarjunaiah.
\newblock An $ hp $-adaptive discontinuous galerkin discretization of a static
  anti-plane shear crack model.
\newblock \emph{arXiv preprint arXiv:2411.00021}, 2024.

\end{thebibliography}

\end{document}